\documentclass[12pt]{article}

\usepackage{amsmath}
\usepackage{amssymb}
\usepackage{amsthm}
\usepackage{latexsym}
\usepackage{color}
\usepackage{graphicx}
\usepackage{color}
\usepackage{caption}
\usepackage{subcaption}
\usepackage{hyperref}
\usepackage{enumerate}
\usepackage{mathtools}
\usepackage{accents}
\usepackage{ulem}
\usepackage{mathrsfs}
\usepackage{authblk}

\DeclareSymbolFont{calletters}{OMS}{cmsy}{m}{n}
\DeclareSymbolFontAlphabet{\mathcal}{calletters}

\newtheorem{Theorem}{Theorem}[section]
\newtheorem{Definition}[Theorem]{Definition}
\newtheorem{Proposition}[Theorem]{Proposition}

\newtheorem{Assumption}[Theorem]{Assumption}

\newtheorem{Lemma}[Theorem]{Lemma}
\newtheorem{Corollary}[Theorem]{Corollary}
\newtheorem{Remark}[Theorem]{Remark}

\makeatletter \@addtoreset{equation}{section}

\newcommand{\brak}[1]{\left(#1\right)}    
\newcommand{\crl}[1]{\left\{#1\right\}}   
\newcommand{\edg}[1]{\left[#1\right]}     
\newcommand{\abv}[1]{\left|#1\right|}     
\newcommand{\norm}[1]{\left\|#1\right\|}  
\newcommand{\agb}[1]{\left\langle #1\right\rangle} 

\usepackage{color}

\newcommand{\bm}{\bibitem}
\newcommand{\no}{\noindent}

\def \B{\mathbb{B}}

\def \E{\mathbb{E}}

\def \P{\mathbb{P}}

\def \R{\mathbb{R}}
\def \Sb{\mathbb{S}}

\def \N{\mathbb{N}}

\def\Ac{{\cal A}}
\def\Bc{{\cal B}}

\def\Dc{{\cal D}}
\def\Ec{{\cal E}}

\def\Pc{{\cal P}}

\def\Sc{\mathcal{S}}
\def\Tc{{\cal T}}

\def\Xc{{\cal X}}
\def\Yc{{\cal Y}}

\def \As{\mathscr{A}}
\def \Es{\mathscr{E}}
\def \Hs{\mathscr{H}}

\def\Dch{\widehat{\Dc}}

\def\Lm{\Lambda}

\def\Rt{\widetilde{R}}

\def\Yt{\widetilde{Y}}

\def\Ebh{\widehat{\E}}

\def \a{\alpha}
\def \ah{\frac{\alpha}{2}}

\def \b{\beta}
\def \bh{\frac{\b}{2}}

\def \Om{\Omega}
\def \om{\omega}

\def \eps{\varepsilon}
\def \xb{\mathbf{x}}

\def \0{\mathbf{0}}
\def \1{\mathbf{1}}

\def \Yb{\overline{Y}}
\def \Yh{\widehat{Y}}
\def \Yu{\underline{Y}}
\def \Zh{\widehat{Z}}
\def \Zb{\overline{Z}}
\def \Kb{\overline{K}}

\def \Hb{\overline{H}}

\def \ab{\bar{a}}

\def \Kt{\tilde{K}}

\def \vp{\varphi}
\def \x{\times}
\def \sigmah{\overline{\sigma}}
\def \sigmat{\widetilde{\sigma}}

\def \half{\frac{1}{2}}

\def \lm{\lambda}

\def \ha{\widehat{\a}}
\def \hb{\widehat{\b}}

\def \Ah{\hat{A}}

\def \yb{\bar{y}}

\def \at{\tilde{a}}

\def \dt{\delta_\theta}
\def \dtt{\widetilde{\delta}_\theta}
\def \dta{\abv{\delta}_\theta}

\def \fb{\bar{f}}

\def \Lu{\uline{L}}
\def \hbu{\underline{\hbar}}
\def \Uu{\underline{U}}
\def \Ru{\underline{R}}
\def \Ybu{\underline{\Yb}}

\def \Hbu{\underline{\Hb}}

\def \Dt{\widetilde{D}}
\def \Ds{\mathscr{D}}

\def \Ab{\overline{A}}

\def \Cs{\mathscr{C}}
\def \bx{\overline{x}}

\def \sb{{\rm{\bar{s}}}}
\def \tb{{\rm{\bar{t}}}}

\def \hv{\vec{\hbar}}
\def \lv{\vec{l}}
\def \uv{\vec{u}}

\def \fv{\vec{f}}
\def \cu{\underline{c}}
\def \cb{\overline{c}}
\def \Dsh{\widehat{\Ds}}

\title{$G$-BSDEs with mean constraints in time-dependent intervals}

\author[a]{Zihao Gu}
\author[a]{Hui Zhao\thanks{Corresponding author: Hui Zhao (amyzh9201@sjtu.edu.cn)}}
{
\affil[a]{School of Mathematical Sciences, Shanghai Jiao Tong University, Shanghai, 200240, China}
}
\date{\today}

\begin{document}
	
	\maketitle
	
	\begin{abstract}
		In this paper, we study a collection of mean-reflected backward stochastic differential equations driven by $G$-Brownian motions ($G$-BSDEs), where $G$-expectations are constrained in some time-dependent intervals. To establish well-posedness results, we firstly construct a backward Skorokhod problem with sublinear expectation, and then apply that in the study of doubly mean-reflected $G$-BSDEs involving Lipschitz and quadratic generators under bounded and unbounded terminal conditions. Also we utilize fixed-point argumentations and $\theta$-methods while solving these equations. Finally, we extend the results to multi-dimensional doubly mean-reflected $G$-BSDEs with diagonal generators.

		\vspace{2mm}
		
		\no {\bf Key words.} Backward Skorokhod problem; $G$-expectation;  Doubly mean-reflected $G$-BSDEs.
		
	\end{abstract}
	
	\section{Introduction}

In 1960s, Skorokhod \cite{SPboundary1,SPboundary2} proposed diffusion processes with reflecting boundaries, where a non-decreasing function is added to a stochastic differential equation (SDE) for pushing it upward in a minimal way, i.e. the Skorokhod condition holds. Since then, Skorokhod problems have been researched extensively by many scholars. Chaleyat-Maurel and El Karoui \cite{ChM-ElK} considered reflected equations with continuous functions, while Chaleyat-Maurel, El Karoui and Marchal \cite{discontref} concerned about the discontinuous case. Yamada \cite{Yamada} studied the uniqueness results for one-dimensional reflected SDEs on the real half-line. For multi-dimensional cases, one can refer to Strook and Varadhan \cite{diffusionBryCond}, Tanaka \cite{Tanaka}, Lions and Sznitman \cite{SDErefB} and so on. 
Skorokhod problem with solution lying between two barriers, also named two-sided Skorokhod problem, was proposed by Kruk, Lehoczky, Ramanan, and Shreve \cite{SP0a,SP0aQu}, where the c\`adl\`ag functions are constrained within a constant interval $[0,a]$. Then Skorokhod problems with time-dependent boundaries were developed by Burdzy, Kang and Ramanan \cite{SPtimeinv},   Słomi\'nski and Wojciechowski \cite{SPtimebs}, Slaby \cite{SlabySP,SlabyESP}, etc.. Recently, Li \cite{SP2nonlinearC} put forward a type with two nonlinear constraints.


After the publication of nonlinear backward stochastic differential equation (BSDE) by Pardoux and Peng \cite{PengBSDE}, El Karoui, Kapoudjian, Pardoux, Peng and Quenez \cite{rbsde} introduced relected BSDEs with one barrier, while Cvitanic and Karatzas \cite{rbsdedg} developed those with two barriers. For more excursions into reflected BSDEs, see also Djehiche and Hamad\`ene \cite{mfrbsde}, Essaky and Hassani \cite{general2RBSDE}, Hamad\`ene and Popier \cite{LpBSDE}, Klimsiak \cite{rbsdefiltered}, Peng and Xu  \cite{IncompleteMarket} and the references therein. Considering super-hedging problems,  Briand, Elie and Hu \cite{mrbsdeBEH} initiated mean-reflected BSDEs taking the form:
\begin{equation}\label{eq:MRBSDE}
		Y_t=\xi+\int_t^Tf(s,Y_s,Z_s)ds-\int_t^TZ_s dB_s+R_T-R_t;~
		\E\edg{\ell(t,Y_t)}\geq0,~\forall t\in[0,T].
\end{equation}
Here, $R$ is a deterministic, continuous and non-decreasing compensation process. And the running loss functions $\brak{\ell(t,\cdot)}_{t\in[0,T]}$ are assumed to be strictly increasing and bi-Lipschitz for each $t\in[0,T]$, of which the constraints on expectations simulate a course of risk management. The associated Skorokhod condition is given by $\int_0^T\E\edg{\ell(t,Y_t)}dR_t=0$. Their existence and uniqueness results were obtained under Lipschitz generators and square-integrable terminal values. Afterwards, Hibon, Hu, Lin, Luo, and Wang \cite{MRQBSDE} proposed  mean-reflected BSDEs with quadratic growth. Hu, Moreau and Wang \cite{GeMRBSDE} extended these results to the case where generators depend on marginal distributions of the solution, and terminal values can be assumed unbounded. In addition, doubly mean-reflected SDEs (BSDEs) were developed by Falkowski and Słomiński \cite{MRSDE2C,MRBSDE2C} and Li \cite{DMRbsde}.


Concerning uncertainty issue in finance, Peng \cite{PengMdG,PengGIto,PengNEC} initiated a time-consistent sublinear expectation, called $G$-expectation. Under such framework, $G$-Brownian motion $B$, its quadratic variation process $\agb{B}$ and induced stochastic calculus were established, of which detailed information can be found in Denis, Hu and Peng \cite{DenisHuPeng}, Li, Peng and Song \cite{LPSSuperD}, Li and Peng \cite{LiXPeng}, Song \cite{SongGeva} and so forth. Wherein, Gao \cite{pathGBM} built up the well-posedness of stochastic differential equations driven by $G$-Brownian motion ($G$-SDEs) and Lin \cite{LinGSDE} developed reflected $G$-SDEs. For multi-dimensional case, Lin and Soumanna Hima \cite{GSDEnonconvex} extended the results of \cite{LinGSDE} to where the range of targeted solution is not necessarily convex. Relatively to the classical framework, Hu, Ji, Peng, and Song \cite{HuJiPengSong} solved uniquely  backward stochastic differential equations driven by $G$-Brownian motion ($G$-BSDEs):
\begin{equation}\label{eq:GBSDE}
	Y_t=\xi+\int_t^Tf(s,Y_s,Z_s)ds+\int_t^Tg(s,Y_s,Z_s)d\agb{B}_s-\int_t^TZ_sdB_s-(K_T-K_t),~t\in[0,T],
\end{equation} 
where the generators $f$ and $g$ are uniformly Lipschitz. The solution to $G$-BSDE \eqref{eq:GBSDE} is a triple of processes $(Y,Z,K)$, where $K$ is a non-increasing $G$-martingale starting from 0. Correspondingly, the comparison theorem, Feynman-Kac formula and Girsanov transformation were proved by Hu, Ji, Peng, and Song \cite{ComFKGir}. Multi-dimensional $G$-BSDEs along with a comparison theorem were proved in Liu \cite{LiumdGBSDE} 
Motivated by model uncertainty, Soner, Touzi and Zhang \cite{martG,WP2bsde} deeply cultivated a topic of second order BSDE (2BSDE), which is related to $G$-framework but allows for various representative sets of probability measures.
 Speaking of reflected $G$-BSDEs, Li, Peng and Soumana Hima \cite{1dRGBSDE} studied a group with one lower obstacle, then Li and Peng \cite{UpOb} focused on one upper barrier. In particular, mean-reflected $G$-SDEs (cf. Li \cite{MRGSDE}) and $G$-BSDEs (cf. Liu and Wang \cite{LWMFGB}) with one barrier were established where designated $G$-expectations are restricted. Lately, He and Li \cite{MRGBSDE2C} developed mean-reflected $G$-BSDEs with two nonlinear constraints on $G$-expectations of different running loss functions.  

Especially, Possama\"i and Zhou \cite{2QGBSDE} studied 2BSDEs with quadratic growth. Soonly after, 
 $G$-BSDEs with quadratic growth (Q$G$-BSDEs) were introduced by Hu, Lin and Soumana Hima \cite{QGBSDEbt} with bounded terminal values. Then, Hu, Tang and Wang \cite{QGBSDEubt} consider Q$G$-BSDEs with generators which are additionally either convex or concave in $Z$, and unbounded terminal values. On reflected Q$G$-BSDEs with lower obstacles, one can refer to Cao and Tang \cite{reflectedQGBSDE}. Recently,  Gu, Lin and Xu \cite{GLX} investigated lower-constrained mean-reflected Q$G$-BSDEs in view of the $''$worst-case scenario$''$.

In this paper, we study mean-reflected $G$-BSDEs with model-relevant sublinear expectations restricted in time-dependent intervals, i.e.  doubly mean-reflected $G$-BSDEs (DMR$G$-BSDEs) written as:
\begin{equation*}
	\left\{
    \begin{lgathered}
    	Y_t=\xi+\int_t^Tf(s,Y_s,Z_s)ds-\int_t^TZ_sdB_s-(K_T-K_t)+(R_T-R_t),~t\in[0,T];\\
    	l_t\leq\Ebh\edg{\hbar(t,Y_t)}\leq u_t,~\forall t\in[0,T].
    \end{lgathered}
	\right.
\end{equation*}
Similarly to \eqref{eq:MRBSDE}, $\brak{\hbar(t,\cdot)}_{t\in[0,T]}$ is a family of running loss functions which are strictly increasing and bi-Lipschitz. The two barriers $l$ and $u$ are bounded and continuous on $[0,T]$ such that $\inf_{t\in[0,T]}\brak{u_t-l_t}>0$. As with reflections, $R$ is deterministic, continuous and of bounded variation starting from 0. Then the Skorokhod condition is given by:
\begin{equation*}
	\int_{s}^{t}\brak{\Ebh\edg{\hbar(v,Y_v)}-l_v}dR_v\leq0,~\mbox{and}~~\int_{s}^{t}\brak{\Ebh\edg{\hbar(v,Y_v)}-u_v}dR_v\leq0,~ \forall0\leq s\leq t\leq T.
\end{equation*} 
In the sequel, we denote by DMR$G$-BSDE ($\xi,f,\hbar,l,u)$ the DMR$G$-BSDE with data ($\xi,f,\hbar,l,u)$.

This paper is organized as follows: Section 2 reviews some basics on the (backward) Skorokhod problems with two barriers and the $G$-expectation theory. In Section 3, we develop a backward Skorokhod problem with sublinear expectation. Then we prove well-posedness results for several DMR$G$-BSDEs in Section 4, including a Lipschitz-generator case, and two quadratic-generator cases each under bounded and unbounded terminal values. In section 5, we study related multi-dimensional DMR$G$-BSDEs with diagonal generators.

\section{Preliminaries}

 Let $T>0$ be a fixed time horizon, $0\leq \sb\leq \tb\leq T$ and $\Xc=\R,~\R^n,~n\in\N^*$. We denote the spaces of functions $x:\R^+\to\Xc$ ($x:[\sb,\tb]\to\Xc$ resp.) by:
 
\indent$C(\R^+;\Xc)$ ($C(\sb,\tb;\Xc)$ resp.): the space of continuous functions;

\indent$C_0(\R^+;\Xc)$ ($C_0(\sb,\tb;\Xc)$ resp.): the space of continuous functions starting from the origin;

\indent$C_b(\R^+;\Xc)$ ($C_b(\sb,\tb;\Xc)$ resp.): the space of bounded continuous functions;

\indent$D(\R^+;\Xc)$ ($D(\sb,\tb;\Xc)$ resp.): the space of c\`adl\`ag fuctions;

\indent$BV(\R^+;\Xc)$ ($BV(\sb,\tb;\Xc)$ resp.): the space of c\`adl\`ag functions having bounded variations.
 
   For any $R\in BV(0,T;\R)$, $\abv{R}^{\sb}_{\tb}$ represents the total variation of $R$ on $[\sb,\tb]$, and $\abv{R}_{\tb}$ denotes that on $[0,\tb]$. 
 Also we simplify the notations by $C_B(\sb,\tb):=C_0(\sb,\tb;\R)\cap BV(\sb,\tb;\R)$, $C_B^T:=C_0(0,T;\R)\cap BV(0,T;\R)$ and $C_B^T(\R^n):=C_0(0,T;\R^n)\cap BV(0,T;\R^n)$.

\subsection{Skorokhod problem and backward Skorokhod problem}

\begin{Definition}[\cite{SPtimeinv}. Skorokhod problem with two barriers]\label{definition:SP}
	Assume that $\bx,l,u\in D(\R^+;\R)$ with $l\leq u$ and $l_0\leq \bx_0\leq u_0$. A pair $(x,k)\in D(\R^+;\R^2)$ is called a solution of the Skorokhod problem associated with $\bx$ and barriers $l,u$, i.e. $(x,k)=\Sb\P_l^u(\bx)$ if
	\begin{enumerate}[($i$)]
		\item $x_t=\bx_t+k_t\in[l_t,u_t]$, $\forall t\in\R^+$;
		\item $k_0=0$, for any $0\leq s\leq t$ such that $\inf_{v\in[s,t]}(u_v-l_v)>0$, $k$ is of bounded variation on $[s,t]$ and 
		\begin{equation}\label{eq:SPmin}
			\int_{[s,t]}(x_v-l_v)dk_v\leq0,~\mbox{and}~~\int_{[s,t]}(x_v-u_v)dk_v\leq0.
		\end{equation}
	\end{enumerate}
\end{Definition}
\begin{Remark}\label{remark:SP}
	Condition \eqref{eq:SPmin}, also called the minimality condition, has an equivalent version by Definition 2.2, Proposition 2.3 and Corollary 2.4 in \cite{SPtimeinv}: for each $0\leq s\leq t$,
	\begin{align*}
		k_t\geq k_s,~&\mbox{if}~x_v<u_v ~\mbox{for all}~v\in(s,t];\\
		k_t\leq k_s,~&\mbox{if}~x_v>l_v~~\mbox{for all}~v\in(s,t];
	\end{align*}
	and for every $t\in\R^+$, $\Delta k_t\geq0$ if $x_t<u_t$ and $\Delta k_t\leq0$ if $x_t>l_t$, where $\Delta k_t:=k_t-k_{t-}$.
\end{Remark}

\begin{Theorem}[\cite{SPtimeinv,SPtimebs}. Explicit formula for the Skorokhod problem]
	Suppose that $\bx,l,u\in D(\R^+;\R)$ with $l\leq u$ and $l_0\leq \bx_0\leq u_0$.  Then for each $\bx\in D(\R^+;\R)$, there exists a unique pair $(x,k)\in D(\R^+;\R^2)$ that solves the Skorokhod problem associated with $\bx$ and barriers $l,u$, i.e., $(x,k)=\Sb\P_l^u(\bx)$, such that
	\begin{equation*}
		k_t=-\max\brak{(\bx_0-u_0)^+\wedge\inf_{0\leq v\leq t}(\bx_v-l_v),\sup_{0\leq s\leq t}[(\bx_s-u_s)\wedge\inf_{s\leq v\leq t}(\bx_v-l_v)]},~t\in\R^+.
	\end{equation*}
\end{Theorem}

\begin{Definition}[\cite{MRBSDE2C}. Backward Skorokhod problem with two barriers]\label{definition:BSP}
Assume that $\bx,l,u\in D(0,T;\R)$, $a\in\R$ satisfy that $l\leq u$ and $l_T\leq a\leq u_T$. A pair $(x,k)\in D(0,T;\R^2)$ is called a solution of the backward Skorokhod problem associated with $\bx,a$ and barriers $l,u$ i.e. $(x,k)=\B\Sb\P_l^u(\bx,a)$ if
	\begin{enumerate}[($i$)]
		\item $x_t=a+\bx_T-\bx_t+k_T-k_t\in[l_t,u_t],~t\in[0,T]$;
		\item $k_0=0$, for any $0\leq s\leq t$ such that $\inf_{v\in[s,t]}(u_v-l_v)>0$, $k$ is of bounded variation on $[s,t]$ and
		\begin{equation}\label{eq:BSPmin}
			\int_{[s,t]}(x_{v-}-l_{v-})dk_v\leq0~\mbox{and}~\int_{[s,t]}(x_{v-}-u_{v-})dk_v\leq0.
		\end{equation} 
	\end{enumerate}
\end{Definition}
\begin{Remark}\label{remark:BSP}
	Similarly to Remark \ref{remark:SP}, one can check that the minimality condition \eqref{eq:BSPmin}  is equivalent to: for each $0\leq s\leq t$,
	\begin{align*}
		k_t\geq k_s,~&\mbox{if}~x_{v-}<u_{v-} ~\mbox{for all}~v\in(s,t];\\
		k_t\leq k_s,~&\mbox{if}~x_{v-}>l_{v-}~~\mbox{for all}~v\in(s,t];
	\end{align*}
	and for every $t\in\R^+$, $\Delta k_t\geq0$ if $x_{t-}<u_{t-}$ and $\Delta k_t\leq0$ if $x_{t-}>l_{t-}$. Notice that, for any $t\in[0,T]$,
	$
		x_{t-}=\max\edg{\min(x_t+\Delta \bx_t,u_{t-}),l_{t-}}$. 
\end{Remark}

\begin{Theorem}[\cite{MRBSDE2C}. Explicit formula for the backward Skorokhod problem]\label{theorem:BSPk}
	For all $\bx,l,u\in D(0,T;\R)$, $a\in\R$ such that $l\leq u$ and $l_T\leq a\leq u_T$, there exists a unique solution to the backward Skorokhod problem, namely $(x,k)=\B\Sb\P_l^u(\bx,a)$. Moreover, for any $t\in[0,T)$, 
	\begin{align*}
		k_T-k_t=&-\max\Big\{(\at-u_T)^+\wedge\inf_{t\leq v\leq T}(\at+\bx_{T-}-\bx_{v-}-l_{v-}),\\
		&\qquad\quad\sup_{t\leq s\leq T}[(\at+\bx_{T-}-\bx_{s-}-u_{s-})\wedge\inf_{t\leq v\leq s}(\at+\bx_{T-}-\bx_{v-}-l_{v-})]\Big\}+\at-(a+\Delta \bx_T),
	\end{align*}
	where $\at=\max(\min(a+\Delta \bx_T,u_{T-}),l_{T-})$.
\end{Theorem}

\subsection{$G$-expectation theory}

Let $\Om_T:=C_0(0,T;\R^d)$ endowed with the supremum norm $\norm{\cdot}:=\sup_{t\in[0,T]}\abv{\cdot}$. Denote $\Om_t:=\crl{\om_{\cdot\wedge t}:~\om\in\Om_T}$ for each $t\in[0,T]$.
Then $\Bc(\Om_T)$ ($\Bc(\Om_t)$ resp.) denotes the Borel $\sigma$-algebra of $\Om_T$ ($\Om_t$ resp.).
 And let $B$ be the canonical process on $\Om_T$ such that $B_t(\om):=\om_t$. Set
 $$L_{ip}(\Om_t):=\crl{\vp(B_{t_1},\cdots,B_{t_n}):~n\in\N,~0\leq t_1<\cdots<t_n\leq t,~\vp\in C_{b,Lip}(\R^{n\x d})},~t\in[0,T],$$ 
where $C_{b,Lip}(\R^{n\x d})$ is the space of bounded Lipschitz functions on $\R^{n\x d}$.

On $(\Om_T,L_{ip}(\Om_T))$, Peng \cite{PengMdG,PengNEC} initiated the $G$-expectation $\Ebh\edg{\cdot}:L_{ip}(\Om_T)\to\R$ and the conditional $G$-expectation $\Ebh_t\edg{\cdot}:L_{ip}(\Om_T)\to L_{ip}(\Om_t),~t\in[0,T]$. We are given a monotonic, sublinear function $G:\Sb_d\to\R$, with $\Sb_d$ being the space of all $d\x d$ symmetric matrices. For any $\vp(B_t),~\vp\in C_{b,Lip}(\R^d)$, we have $\Ebh\edg{\vp(B_t)}=u(t,\0)$, where $u(t,\xb)$ is the viscosity solution to the fully nonlinear PDE:
\begin{equation*}
		 \partial_t u(t,\xb)-G(\partial^2_{\xb\xb}u(t,\xb))=0,~\forall(t,\xb)\in[0,T]\x\R^d;~
		 u(0,\xb)=\vp(\xb).
\end{equation*} 
 $G$ is non-degenerate, i.e. there exist two constants $0<\sigmat^{-1}<\sigmah<\infty$ such that
\begin{equation*}
	\frac{1}{2\sigmat^2}{\rm tr}[A-A']\leq G(A)-G(A')\leq\frac12\sigmah^2{\rm tr}[A-A'],~\mbox{for all}~A\geq A'.
\end{equation*}

 Let $p\geq1$. We denote by $L_G^p(\Om_T)$ ($L_G^p(\Om_t)$ resp.) the completion of $L_{ip}(\Om_T)$ ($L_{ip}(\Om_t)$ resp.) under the norm $(\Ebh\edg{\abv{\cdot}^p})^{\frac 1p}$. Hereby, $\Ebh\edg{\cdot}$ and $\Ebh_t\edg{\cdot},~t\in[0,T]$ can be continuously extended to $L_G^1(\Om_T)$. 
 Indeed, $G$-expectation is an upper expectation.
 \begin{Theorem}[\cite{DenisHuPeng}]
 	There exists a weakly compact set $\Pc$ of probability measures on $(\Om_T,L_G^1(\Om_T))$ such that
 	\begin{equation*}
 		\Ebh\edg{\xi}=\sup_{\P\in\Pc}\E^{\P}\edg{\xi},~\mbox{for all}~\xi\in L_G^1(\Om_T).
 	\end{equation*}
 Then $\Pc$ is called a set that represents $\Ebh$.
 \end{Theorem}
For this $\Pc$, we define a capacity as:
\begin{equation*}
	c(A):=\sup_{\P\in\Pc}\P(A),~A\in\Bc(\Om_T).
\end{equation*}
 A set $A\in\Bc(\Om_T)$ is polar if $c(A)=0$. A property holds $''quasi-surely ~(q.s.)''$ if it holds outside a polar set. We do not distinguish between two random variables or processes $X$ and $Y$ if $X=Y~q.s.$.
 
 There are also abundant constructions and theorems under $G$-framework. 
 $L_G^\infty(\Om_T)$ ($L_G^\infty(\Om_t)$ resp.) is the completion of $L_{ip}(\Om_T)$ ($L_{ip}(\Om_t)$ resp.) under the norm $$\norm{\xi}_{L^\infty_G}:=\inf\crl{M\geq0:\abv{\xi}\leq M~q.s.}.$$ 
 Let $0\leq \sb\leq \tb\leq T$. Set 
 \begin{equation*}
 	M_G^0(\sb,\tb):=\crl{\zeta_\cdot=\sum_{j=0}^{n-1}\xi_j\1_{[t_j,t_{j+1})}(\cdot):~n\in\N,~ \sb=t_0<t_1<\cdots<t_n=\tb,~\xi_j\in L_{ip}(\Om_{t_j})},
 \end{equation*}
 and
 \begin{equation*}
 	S_G^0(\sb,\tb):=\crl{h(\cdot,B_{t_1\wedge \cdot},\cdots,B_{t_n\wedge \cdot}):~n\in\N,~ \sb\leq t_1<t_2<\cdots<t_n\leq \tb,~h\in C_{b,Lip}(\R^{1+n\x d})}.
 \end{equation*}
Denote by $M_G^p(\sb,\tb)$ the completion of $M_G^0(\sb,\tb)$ under the norm $(\Ebh[\int_{\sb}^{\tb}\abv{\cdot}^pdv])^{\frac1p}$, and by $H_G^p(\sb,\tb)$ that under $(\Ebh[(\int_{\sb}^{\tb}\abv{\cdot}^2dv)^{\frac p2}])^{\frac1p}$. Also, the completion of $S_G^0(\sb,\tb)$ under the norm $(\Ebh[\sup_{v\in[\sb,\tb]}\abv{\cdot}^p])^{\frac 1p}$ is symbolized by $S_G^p(\sb,\tb)$. Moreover, we denote by $S_G^{\infty}(\sb,\tb)$ the completion of $S_G^0(\sb,\tb)$ under the norm
 $\norm{\eta}_{S_G^{\infty}(\sb,\tb)}:=\norm{\sup_{v\in[\sb,\tb]}\abv{\eta_v}}_{L_G^\infty}.$
 
 On such basis, $\Ec_G^p(\sb,\tb)$ is defined as the collection of stochastic processes $Y$ such that $e^{Y}\in S_G^p(\sb,\tb)$, along with
$\Es_G^T:=\bigcap_{p\geq1}\Ec_G^p(0,T)$. Then 
 $\Ac_G^p(\sb,\tb)$ stands for the space of all non-increasing $G$-martingales $K$ such that $K_{\sb}=0$ and $K_{\tb}\in L_G^p(\Om_{\tb})$. Correspondingly, 
$\Hs_G^T:=\bigcap_{p\geq1}H_G^p(0,T)$, and $\As_G^T:=\bigcap_{p\geq1}\Ac_G^p(0,T)$. The notation $C_G^p(\sb,\tb)$ will be utilized frequently to symbolize the space of processes $Y$ such that $Y_v\in L_G^p(\Om_v),~\forall v\in[\sb,\tb]$ and the mapping $v\mapsto Y_v$ is continuous under the norm $\norm{\cdot}_{L_G^1}$.
Besides, we put $\Sc_G^p(\sb,\tb):=S_G^p(\sb,\tb)\x H_G^p(\sb,\tb;\R^d)\x \Ac_G^p(\sb,\tb)$, for $p\geq1$.

 It is worth mentioning that, the canonical process $B_\cdot:=\brak{B^i_\cdot}_{1\leq i\leq d}$ is a $d$-dimensional $G$-Brownian motion on the $G$-expectation space $(\Om_T,L_G^1(\Om_T),\Ebh[\cdot],(\Ebh_t[\cdot])_{t\in[0,T]})$. For each $1\leq i,j\leq d$, we denote by $\agb{B^i,B^j}$ the mutual variation process. Following Corollary 3.5.8 of Peng \cite{PengNEC}, it holds that
 \begin{equation*}
 	\sigmat^{-2}I_d\leq \edg{\frac{d\agb{B^i,B^j}_t}{dt}}_{i,j=1}^d\leq\sigmah^2I_d.
 \end{equation*}
Then given two processes $\eta\in M_G^1(\sb,\tb)$ and $\zeta\in H_G^1(\sb,\tb)$, $G$-It\^o integrals $\int\eta_vd\agb{B^i,B^j}_v$ and $\int\zeta_vdB^i_v$ could be constructed by a standard approximation procedure. Regarding the multi-dimensional case, $M_G^p(\sb,\tb;\R^n)$ refers to the space containing all $\R^n$-valued processes of which every component belongs to $M_G^p(\sb,\tb)$. Similarly, we can define $L_G^p(\Om_t;\R^n)$, $H_G^p(\sb,\tb;\R^d)$, $S_G^p(\sb,\tb;\R^n)$, $\Ec_G^p(\sb,\tb;\R^n)$, $\Ac_G^p(\sb,\tb;\R^n)$, $S_G^\infty(\sb,\tb;\R^n)$, and $$\Sc_G^p(\sb,\tb;\R^n):=S_G^p(\sb,\tb;\R^n)\x H_G^p(\sb,\tb;\R^{n\x d})\x \Ac_G^p(\sb,\tb;\R^n).$$ Then we denote $\Hs_G^T(\R^d):=\bigcap_{p\geq1}H_G^p(0,T;\R^d)$, thus $\Es_G^T(\R^n)$, $\Hs_G^T(\R^{n\x d})$, and $\As_G^T(\R^n)$ are defined respectively. For simplicity, we set
\begin{equation*}
	\int_{\sb}^{\tb}\eta_vd\agb{B}_v:=\sum_{i,j=1}^d\int_{\sb}^{\tb}\eta_v^{ij}d\agb{B^i,B^j}_v,~\mbox{and}~~\int_{\sb}^{\tb}\zeta_vdB_v:=\sum_{i=1}^d\int_{\sb}^{\tb}\zeta^i_vdB^i_v,
\end{equation*}
for any $\eta\in M_G^1(\sb,\tb;\R^{d\x d})$ and $\zeta\in H_G^1(\sb,\tb;\R^d)$.

\begin{Theorem}[\cite{DenisHuPeng,QGBSDEubt}. Monotone convergence theorem and Fatou's lemma]\label{theorem:mcvg}
	Let $X,~X_n,~n\geq1$ be $\Bc(\Om_T)$-measurable random variables.
	\begin{enumerate}[($i$)]
		\item If $\crl{X_n}_{n=1}^\infty\subset L_G^1(\Om_T)$ satisfies that $X_n\downarrow X$ q.s., then $\Ebh\edg{X_n}\downarrow\Ebh\edg{X}$.
		\item If $X_n\geq 0,~\forall n\geq1$. Then, $\Ebh\edg{\liminf_{n\to\infty}X_n}\leq \liminf_{n\to\infty}\Ebh\edg{X_n}$.
\end{enumerate} 
\end{Theorem}

\begin{Lemma}[\cite{LiumdGBSDE}. Doob's maximal inequality for $G$-martingales]
	For any $\a\geq 1$, $\delta>0$, and for all $\xi\in L_G^{\a+\delta}(\Om_T)$, there exists some constant $C(\a,\delta)>0$ depending only on $\a$ and $\delta$ such that
	\begin{equation*}
		\Ebh\edg{\sup_{t\in[0,T]}\Ebh_t\edg{\abv{\xi}^\a}}\leq C(\a,\delta)\brak{\Ebh\edg{\abv{\xi}^{\a+\delta}}}^{\frac{\a}{\a+\delta}}.
	\end{equation*}
\end{Lemma}
\begin{Corollary}[\cite{GLX,QGBSDEubt}. Application of Doob's maximal inequality]\label{corollary:Doob}
	If $X$ is positive and $e^X\in L_G^2(\Om_T)$, then there exists a constant $\Ah_G>0$ such that (taking $\a=\delta=2$)
	\begin{equation*}
		\Ebh\edg{\sup_{t\in[0,T]}\Ebh_t\edg{e^X}}=\Ebh\edg{\sup_{t\in[0,T]}\Ebh_t\edg{\abv{e^{\frac X2}}^2}}\leq \Ah_G\Ebh\edg{e^{2X}}.
	\end{equation*}
\end{Corollary}

\section{Backward Skorokhod problem with sublinear expectation}
In this section, we define a backward Skorokhod problem with sublinear expectation, derive the corresponding explicit formula and deduce some useful a priori estimates. 
Throughout this paper, we fix the following assumptions on reflections.
\begin{Assumption}\label{assumption:HR}
	 The obstacles $l,u$ and the running loss function $\hbar:[0,T]\x\Om_T\x\R\to\R$ satisfy:
	\begin{enumerate}[($i$)]
		\item $l,u\in C_b(0,T;R)$ with $\inf_{t\in[0,T]}\brak{u_t-l_t}>0.$
		\item $(t,y)\mapsto\hbar(t,y)$ is uniformly continuous, uniformly in $\om\in\Om$.
		\item For any $t\in[0,T]$, $y\mapsto\hbar(t,y)$ is strictly increasing.
		\item For any $t\in[0,T]$, $\hbar(t,y)\in L_G^1(\Om_T)$, $\forall y\in\R$,
		$$\Ebh\edg{\lim_{y\downarrow-\infty}\hbar(t,y)}<\inf_{s\in[0,T]}l_s<\sup_{s\in[0,T]}u_s< \Ebh\edg{\lim_{y\uparrow\infty}\hbar(t,y)}.$$
		\item For any $t\in[0,T]$, there exists some $C_\hbar>0$ such that $$\abv{\hbar(t,y)}\leq C_\hbar(1+\abv{y}),~\forall y\in\R.$$
		\item For any $t\in[0,T]$, there exist some constants $0<\cu\leq \cb$ such that for any $(y,y')\in\R^2$,
		$$\cu\abv{y-y'}\leq \abv{\hbar(t,y)-\hbar(t,y')}\leq \cb\abv{y-y'}.$$
	\end{enumerate}
\end{Assumption}

For $t\in[0,T]$, let $X\in L_G^1(\Om_t)$ and $\hbar$ satisfies ($H_r$), we define a mapping $\Hb(t,\cdot,X):\R\to\R$:
\begin{equation*}
	\Hb(t,x,X)=\Ebh\edg{\hbar(t,x+X-\Ebh\edg{X})},~x\in\R.
\end{equation*}
By Assumption ($H_r$)-($iv$)($v$), $\Hb(t,\cdot,X)$ is well-defined.  

\begin{Proposition}\label{proposition:hbarpropoty}
	Let Assumption \ref{assumption:HR} hold, $X\in L_G^1(\Om)$, $\zeta\in M_G^1(0,T)$, and $S\in S_G^1(0,T)$. 
	\begin{enumerate}[($i$)]
		\item For each $(t,x)\in[0,T]\x\R$, $\hbar(t,x+X-\Ebh\edg{X})\in L_G^1(\Om_T)$.
		\item The mapping $x\mapsto\hbar(t,x+X-\Ebh\edg{X})$ is continuous. Moreover, $x\mapsto\Hb(t,x,X)$ is continuous, strictly increasing and
		$$\lim_{x\downarrow-\infty}\Hb(t,x,X)=-\infty,~ \lim_{x\uparrow\infty}\Hb(t,x,X)=+\infty.$$
		\item The mapping $t\mapsto\hbar\brak{t,x+\Ebh_t\edg{X}-\int_0^t\zeta_sds-\Ebh\edg{X-\int_0^t\zeta_sds}}$ is continuous under the norm $\norm{\cdot}_{L_G^1}$. Moreover, $t\mapsto\Hb(t,x,\Ebh_t\edg{X}-\int_0^t\zeta_sds)$ is continuous.
		\item The mapping $t\mapsto\hbar(t,x+S_t-\Ebh\edg{S_t})$ is continuous under the norm $\norm{\cdot}_{L_G^1}$. Moreover, $t\mapsto\Hb(t,x,S_t)$ is continuous. 
		\item For any $t\in[0,T]$, the mapping $\Hb(t,\cdot,X)$ is bi-Lipschitz with $\cu,\cb>0$, i.e. 
		\begin{equation*}
			\cu\abv{x-x'}\leq\abv{\Hb_(t,x,X)-\Hb(t,x',X)}\leq \cb\abv{x-x'},~x,x'\in\R.
		\end{equation*}
	\end{enumerate}
\end{Proposition}
\begin{Remark}
Proposition \ref{proposition:hbarpropoty} can be proved from Lemma 3.3 in \cite{MRGBSDE2C} and Lemma 3.9 in \cite{MRGSDE}, where the assumptions on running loss functions coincide with Assumption \ref{assumption:HR}-($ii$)($iii$)($v$)($vi$). Note that, by Proposition 3.3 in \cite{LWMFGB} and Lemma 3.7 in \cite{LPSSuperD}, assertions ($i$)-($iii$)-($iv$) can be synthesized as: 
\begin{itemize}
	\item If $S\in C_G^1(0,T)$, i.e., the mapping $t\mapsto S_t$ is continuous under the norm $\norm{\cdot}_{L_G^1}$ and $S_t\in L_G^1(\Om_t),~\forall t\in[0,T]$, so does $\hbar(\cdot,x+S_\cdot-\Ebh\edg{S_\cdot})$.
\end{itemize}
\end{Remark}

\begin{Definition}\label{definition:BSPSE}
	Let Assumption \ref{assumption:HR} hold, and $\Yb\in C_G^1(0,T)$ with $l_T\leq\Ebh\edg{\hbar(T,\Yb_T)}\leq u_T$. A couple of processes $(Y,R)\in C_G^1(0,T)\x C_B^T$ is called a solution to the backward Skorokhod problem with sublinear expectation associated to $\Yb,\hbar$ and two barriers $l,u$,  i.e. $(Y,R)=\widehat{\B\Sb\P}_l^u(\Yb,\hbar)$, if
	\begin{enumerate}[($i$)]
		\item $Y_t=\Yb_t+R_T-R_t,~t\in[0,T]$;
		\item $l_t\leq\Ebh\edg{\hbar(t,Y_t)}\leq u_t,~\forall t\in[0,T]$;
		\item $R_0=0$, $R$ is deterministic, of bounded variation and for any $0\leq s\leq t\leq T$, it holds that
		\begin{equation}\label{eq:minBSPE}
			\int_{s}^{t}\brak{\Ebh\edg{\hbar(v,Y_v)}-l_v}dR_v\leq0,~\mbox{and}~~\int_{s}^{t}\brak{\Ebh\edg{\hbar(v,Y_v)}-u_v}dR_v\leq0.
		\end{equation}
	\end{enumerate}
\end{Definition}
\begin{Lemma}\label{lemma:variii}
	The minimality condition \eqref{eq:minBSPE} is equivalent to: 
	\begin{align}
		\mbox{for each}~~0\leq s<t\leq T,&\label{eq:iii}\\ 
		&R_{t}\geq R_{s},~\mbox{if}~\Ebh\edg{\hbar(v,Y_v)}<u_v~\mbox{for all}~v\in(s,t];\notag\\
		&R_{t}\leq R_{s},~\mbox{if}~\Ebh\edg{\hbar(v,Y_v)}>l_v~\mbox{for all}~v\in(s,t];\notag\\
		\mbox{and for any}~~ v\in[0,T], ~&dR_v\geq0~\mbox{ if}~~\Ebh\edg{\hbar(v,Y_v)}<u_v~\mbox{and}~~dR_v\leq0~\mbox{ if}~~\Ebh\edg{\hbar(v,Y_v)}>l_v.\notag 
	\end{align}
\end{Lemma}
\noindent{\bf\underline{Proof.}} The property $''$\eqref{eq:minBSPE}$\Longrightarrow$\eqref{eq:iii}$''$ can be easily obtained. So we will focus on the other direction $''$\eqref{eq:iii}$\Longrightarrow$\eqref{eq:minBSPE}$''$ following Proposition 2.3 in \cite{SPtimeinv}. Suppose that $(Y,R)\in C_G^1(0,T)\x C_B^T$ satisfy conditions ($i$)($ii$) with condition \eqref{eq:iii}, $R_0=0$, and $R$ is of bounded variation. Then the Lebesgue-Stieltjes measure $dR$ is absolutely continuous w.r.t. the corresponding total variation measure $d\abv{R}$. Denote the associated Radon-Nikodym derivative by $g:=\frac{dR}{d\abv{R}}$. Then $g$ is $d\abv{R}$-measurable, $g(s)\in\crl{-1,1}$ for $d\abv{R}-a.e.~ s\in[0,T]$ and
$	R_t=\int_0^tg(s)d\abv{R}_s$.
Moreover, for $d\abv{R}-a.e.~ s\in[0,T)$,
\begin{equation}\label{eq:RNderi}
	g(s)=\lim_{n\to\infty}\frac{R_{s+\eps_n}-R_s}{\abv{R}_{s+\eps_n}-\abv{R}_s},
\end{equation}
where $\crl{\eps_n}_{n\geq1}\subset\R^+$ depending on $s$ is such that $\abv{R}_{s+\eps_n}>\abv{R}_s$ and $\lim_{n\to\infty}\eps_n=0$. Then for each $t\in[0,T]$, we can define
\begin{equation*}
	R^l_t:=\int_0^t\1_{\crl{g(s)=1}}d\abv{R}_s~\mbox{and}~R^u_t:=\int_0^t\1_{\crl{g(s)=-1}}d\abv{R}_s,
\end{equation*}
Obviously, $R=R^l-R^u$ and $R^l,R^u$ are non-decreasing with $R^l_0=R^u_0=0$. It follows that
\begin{equation*}
	\int_0^T\1_{\crl{\Ebh\edg{\hbar(s,Y_s)}>l_s}}dR^l_s=\int_0^T\1_{\crl{\Ebh\edg{\hbar(s,Y_s)}>l_s}}\1_{\crl{g(s)=1}}d\abv{R}_s\geq0.
\end{equation*}
Assume that there exists some $s\in[0,T]$ such that $\Ebh\edg{\hbar(s,Y_s)}>l_s$, $g(s)=1$ and \eqref{eq:RNderi} is well-defined. From the continuity of $t\mapsto\Ebh\edg{\hbar(t,Y_t)}$ and $l$, there exists some $\delta_s>0$ such that $\Ebh\edg{\hbar(v,Y_v)}>l_v$ for all $v\in[s,s+\delta_s]$. By Assertion \eqref{eq:iii}, we have $R_v\leq R_s,~\forall v\in[s,s+\delta_s]$. However, since $g(s)=1$ and $\abv{R}$ is non-decreasing, for some $N_s\in\N$ sufficiently large, it holds that $\eps_n<\delta_s$ and $R_{s+\eps_n}>R_s,~\forall n\geq N_s$. Hence a contradiction occurs. In other words, either $\crl{\Ebh\edg{\hbar(s,Y_s)}>l_s}$ or $\crl{g(s)=1}$ fails for $d\abv{R}-a.e. ~s\in[0,T]$. Therefore for each $0\leq s\leq t\leq T$,
\begin{align*}
	\int_{s}^{t}\brak{\Ebh\edg{\hbar(v,Y_v)}-l_v}dR_v&=\int_{s}^{t}\brak{\Ebh\edg{\hbar(v,Y_v)}-l_v}dR_v^l-\int_{s}^{t}\brak{\Ebh\edg{\hbar(v,Y_v)}-l_v}dR_v^u\\
	&=-\int_{s}^{t}\brak{\Ebh\edg{\hbar(v,Y_v)}-l_v}dR_v^u\leq0.
\end{align*}
Hereby, the second inequality in condition \eqref{eq:minBSPE} can be derived analogously.  
\qed

For each $t\in[0,T]$, we define the operators $L_t:L_G^1(\Om_t)\to\R$ and $U_t:L_G^1(\Om_t)\to\R$:
\begin{equation*}
	L_t(X):=\Hb{}^{-1}(t,\cdot,X)(l_t),
~\mbox{and}~~
	U_t(X):=\Hb{}^{-1}(t,\cdot,X)(u_t),~X\in L_G^1(\Om_t).
\end{equation*}
Obviously, Assumption \ref{assumption:HR}-($i$) and Proposition \ref{proposition:hbarpropoty}-($ii$) imply that $L_t(\cdot)$ and $U_t(\cdot)$ are well-defined. So we are in a position to study the related properties:
 
\begin{Proposition}\label{proposition:ULcon}
	Let Assumption \ref{assumption:HR} hold. 
	\begin{enumerate}[($i$)]
		\item If $S\in C_G^1(0,T)$, then the mappings $t\mapsto L_t(S_t)$ and $t\mapsto U_t(S_t)$ are continuous. In particular, the mappings $t\mapsto L_t(0)$ and $t\mapsto U_t(0)$ are continuous and bounded.
		\item If $X,X'\in L_G^1(\Om_t)$, then
		\begin{equation*}
			\abv{L_t(X)-L_t(X')}\leq\frac{2\cb}{\cu}\Ebh\edg{\abv{X-X'}},
		~\mbox{and}~~
	    	\abv{U_t(X)-U_t(X')}\leq\frac{2\cb}{\cu}\Ebh\edg{\abv{X-X'}}.
	    \end{equation*} 
	\end{enumerate}	
\end{Proposition}
\noindent{\bf\underline{Proof.}}
We only prove the assertions on $L_\cdot(\cdot)$, as those on $U_\cdot(\cdot)$ can be obtained similarly.

\noindent($i$) Fix $t\in[0,T]$ and $x,x'\in\R$ arbitrarily chosen such that $x<L_t(S_t)<x'$. 
Set $\Lm:=\min\brak{\abv{l_t-\Hb(t,x,S_t)},\abv{l_t-\Hb(t,x',S_t)}}>0$. By Assumption \ref{assumption:HR}-($i$) and Proposition \ref{proposition:hbarpropoty}-($iv$), there exists some $\delta_\Lm>0$ such that for any $s\in[0,T]\cap[t-\delta_\Lm,t+\delta_\Lm]$,
\begin{align*}
	\abv{l_t-l_s}\leq\frac{\Lm}{3},~\abv{\Hb(s,x,S_s)-\Hb(t,x,S_t)}\leq\frac{\Lm}{3},~\mbox{and}~~\abv{\Hb(s,x',S_s)-\Hb(t,x',S_t)}\leq \frac{\Lm}{3},
\end{align*}
thus
\begin{align*}
	\Hb(s,x,S_s)&\leq\Hb(t,x,S_t)+\frac{\Lm}{3}\leq l_t-\frac{2\Lm}{3}\leq l_s-\frac{\Lm}{3}<l_s\leq l_t+\frac{\Lm}{3}\leq \Hb(t,x',S_t)-\frac{2\Lm}{3}\\
	&\leq \Hb(s,x',S_s)-\frac{\Lm}{3}<\Hb(s,x',S_s),
\end{align*}
which implies that $L_s(S_s)\in(x,x')$.

\noindent($ii$) For any $t\in[0,T]$, by Proposition \ref{proposition:hbarpropoty}-($v$), we have
\begin{align*}
	&\quad\Hb\brak{t,L_t(X')+\frac{2\cb}{\cu}\Ebh\edg{\abv{X-X'}},X}=\Ebh\edg{\hbar(t,L_t(X')+\frac{2\cb}{\cu}\Ebh\edg{\abv{X-X'}}+X-\Ebh\edg{X})}\\
	&\geq 2\cb\Ebh\edg{\abv{X-X'}}+\Ebh\edg{\hbar(t,L_t(X')+X-\Ebh\edg{X})}\\
	&\geq 2\cb\Ebh\edg{\abv{X-X'}}+\Ebh\edg{\hbar(t,L_t(X')+X'-\Ebh\edg{X'})}-\cb\Ebh\edg{\abv{X-X'}}-\cb\abv{\Ebh\edg{X}-\Ebh\edg{X'}}\\
	&\geq \Hb(t,L_t(X'),X')=l_t.
\end{align*} 
Then by Proposition \ref{proposition:hbarpropoty}
-($ii$), we have
\begin{equation*}
	L_t(X)\leq L_t(X')+\frac{2\cb}{\cu}\Ebh\edg{\abv{X-X'}}.
\end{equation*}
By symmetry we get the desired result.
\qed

\begin{Theorem}\label{theorem:BSPE&U}
	Let Assumption \ref{assumption:HR} hold. If  $\Yb\in C_G^1(0,T)$ with $l_T\leq\Ebh\edg{\hbar(T,\Yb_T)}\leq u_T$, then there exists a unique solution to the backward Skorokhod problem with sublinear expectation, i.e. $C_G^1(0,T)\x C_B^T\ni(Y,R)=\widehat{\B\Sb\P}_l^u(\Yb,\hbar)$. Moreover, for $t\in[0,T]$,
	\begin{align*}
		R_T-R_t&=-\max\bigg\{(\Ebh\edg{\Yb_T}-U_T(\Yb_T))^+\wedge\inf_{t\leq v\leq T}(\Ebh\edg{\Yb_v}-L_v(\Yb_v)),\\
		&~~~~~~~~~~~~~~~~~~~~~~~~~~~~~~~~~~~~~~~~~\sup_{t\leq s\leq T}\Big[(\Ebh\edg{\Yb_s}-U_s(\Yb_s))\wedge\inf_{t\leq v\leq s}(\Ebh\edg{\Yb_v}-L_v(\Yb_v))\Big]\bigg\}.
	\end{align*}
\end{Theorem}

\noindent{\bf\underline{Proof.}} 
Set $\ab=\Ebh\edg{\Yb_T}$, $\yb_T=-\Ebh\edg{\Yb_T}$, $\yb_t=-\Ebh\edg{\Yb_t}$. Clearly, we have $\yb\in C(0,T;\R)$. 
Note that 
$$\Hb(T,\ab,\Yb_T)=\Ebh\edg{\hbar(T,\Yb_T-\Ebh\edg{\Yb_T}+\ab)}=\Ebh\edg{\hbar(T,\Yb_T)}\in[l_T,u_T],$$
along with the fact that $\Hb(T,\cdot,\Yb_T)$ is strictly increasing, shows that
$$\ab\in[L_T(\Yb_T),U_T(\Yb_T)].$$
From Assumption \ref{assumption:HR}-($i$) and strict increasing property of $\Hb(t,\cdot,\Yb_t),~\forall t\in[0,T]$, we have
\begin{equation}\label{eq:infUL}
	\inf_{t\in[0,T]}\brak{U_t(\Yb_t)-L_t(\Yb_t)}>0.
\end{equation}
Applying Theorem \ref{theorem:BSPk}, there exists a unique solution to the deterministic backward Skorokhod problem $(y,R)=\B\Sb\P_{L_\cdot(\Yb_\cdot)}^{U_\cdot(\Yb_\cdot)}(\yb,\ab)$ such that 
\begin{equation*}
	y_t=\ab+\yb_T-\yb_t+R_T-R_t\in[L_t(\Yb_t),U_t(\Yb_t)],~\forall t\in[0,T],
\end{equation*}
$R_0=0$, and
\begin{align*}
	&R_T-R_t=-\max\bigg\{(\ab-U_T(\Yb_T))^+\wedge\inf_{t\leq v\leq T}(\ab+\yb_T-\yb_v-L_v(\Yb_v)),\\
	&~~~~~~~~~~~~~~~~~~~~~~~~~~~~~~~~~~~~~~~~~~\sup_{t\leq s\leq T}\Big[(\ab+\yb_T-\yb_s-U_s(\Yb_s))\wedge\inf_{t\leq v\leq s}(\ab+\yb_T-\yb_v-L_v(\Yb_v))\Big]\bigg\},
\end{align*}
thus the explicit formula of $R$ holds.
Therefore, $t\mapsto y_t$ is continuous. Letting
\begin{equation*}
	Y_t=\Yb_t+R_T-R_t,~t\in[0,T],
\end{equation*}
and one can check that $y_t=\Ebh\edg{Y_t}$ and $(Y,R)=\widehat{\B\Sb\P}_l^u(\Yb,\hbar)$. Actually, since $\Hb(t,\cdot,\Yb_t)$ is strictly increasing, for any $t\in[0,T]$,
\begin{align*}
	\Ebh\edg{\hbar(t,Y_t)}=\Ebh\edg{\hbar(t,\Yb_t+R_T-R_t)}=\Hb(t,y_t,\Yb_t)\in[\Hb(t,L_t(\Yb_t),\Yb_t),\Hb(t,U_t(\Yb_t),\Yb_t)]=[l_t,u_t],
\end{align*}
then condition ($ii$) of Definition \ref{definition:BSPSE} holds.
For condition ($iii$), it is sufficient to note that
\begin{equation}\label{eq:equiv}
	\left\{\begin{lgathered}
		y_t=\Ebh\edg{Y_t}<U_t(\Yb_t)~\Longleftrightarrow~\Ebh\edg{\hbar(t,Y_t)}<u_t,\\
		y_t=\Ebh\edg{Y_t}>L_t(\Yb_t)~\Longleftrightarrow~\Ebh\edg{\hbar(t,Y_t)}>l_t.
	\end{lgathered}
	\right.
\end{equation}
Referring to  Remark \ref{remark:BSP}, we can deduce that condition \eqref{eq:iii} holds, i.e. condition ($iii$) holds by Lemma \ref{lemma:variii}. 

Next, we prove uniqueness. Suppose that there exists another solution $(Y',R')=\widehat{\B\Sb\P}_l^u(\Yb,\hbar)$ and $Y'_t=\Yb_t+R'_T-R'_t,~t\in[0,T]$.
There exists some $\tb_0\in[0,T]$ such that $Y_{\tb_0}\neq Y'_{\tb_0}$, and we can say $Y_{\tb_0}>Y'_{\tb_0}~ q.s.$ in that $R$ and $R'$ are deterministic. Thus $R_T-R_{\tb_0}>R'_T-R'_{\tb_0}$. Define 
$$\tb_1:=\inf\crl{s\geq \tb_0:~R_T-R_s=R'_T-R'_s}\wedge T.$$
 By the 
 fact that $x\mapsto\Ebh\edg{\hbar(s,\Yb_s+x)}$ is strictly increasing which is implied by Proposition 3.3 in \cite{LWMFGB}, 
 we have $\Ebh\edg{\hbar(s,Y'_s)}<\Ebh\edg{\hbar(s,Y_s)}\leq u_s$, and $dR'_s\geq0$, for all $s\in(\tb_0,\tb_1)$, in view of assertion \eqref{eq:iii}. It follows that,
\begin{equation}\label{eq:Rcomp}
	R_T-R_{\tb_1}=R'_T-R'_{\tb_1}\leq R'_T-R'_{\tb_0}<R_T-R_{\tb_0}.
\end{equation}
Meanwhile, for all $s\in(\tb_0,\tb_1)$, we have $l_s\leq \Ebh\edg{\hbar(s,Y'_s)}<\Ebh\edg{\hbar(s,Y_s)}$, therefore $dR_s\leq0$ and
\begin{equation*}
	R_T-R_{\tb_1}\geq R_T-R_{\tb_0},
\end{equation*}
contradicting \eqref{eq:Rcomp}.
\qed

\begin{Remark}\label{remark: Rflat}
	Let $\Yb,\hbar,l,u$ be in Definition \ref{definition:BSPSE} and $(Y,R)=\widehat{\B\Sb\P}_l^u(\Yb,\hbar)$. For each $t\in[0,T]$ such that $l_t<\Ebh\edg{\hbar(t,Y_t)}<u_t$, i.e. $L_t(\Yb_t)<\Ebh\edg{Y_t}<U_t(\Yb_t)$, we have $dR_t=0$.
	
	Moreover, if there exist $0\leq s<t\leq T$ such that $l_v<\Ebh\edg{\hbar(v,Y_v)}<u_v$, i.e. $L_v(\Yb_v)<\Ebh\edg{Y_v}<U_v(\Yb_v),~\forall v\in[s,t]$, then $R_s=R_t$.
\end{Remark}

\begin{Remark}
	Note that an alternative construction of $R$ can be started by putting
	$$\Hb'(s,x,X):=\Ebh\edg{\hbar(s,X+\Ebh\edg{-X}+x)},~x\in\R,~X\in L_G^1(\Om_s),~s\in[0,T].$$
	Then Proposition \ref{proposition:hbarpropoty} holds for $\Hb'$ and the corresponding operators $L'_s(\cdot)$ and $U'_s(\cdot)$ can be defined satisfying Proposition \ref{proposition:ULcon}. However, the proof for uniqueness is irrelevant to selections for $(\ab,\yb,L_\cdot(\Yb_\cdot),U_\cdot(\Yb_\cdot))$. In the following, we derive a priori estimates via equivalence \eqref{eq:equiv} only. 
\end{Remark}

\begin{Proposition}\label{proposition:BSPSEaprDf}
	For $i=1,2$, let $\Yb^i\in C_G^1(0,T)$ with $l_T\leq\Ebh\edg{\hbar^i(T,\Yb^i_T)}\leq u_T$, and Assumption \ref{assumption:HR} hold for $l^i,u^i$ and $\hbar^i$ with the same positive constants $C_\hbar,\cu,\cb$. If $C_G^1(0,T)\x C_B^T\ni(Y^i,R^i)=\widehat{\B\Sb\P}_{l^i}^{u^i}(\Yb^i,\hbar^i)$, then for any $t\in[0,T]$,
	\begin{align}
		&\quad\abv{(R^1_T-R^1_t)-(R^2_T-R^2_t)}\notag\\
		&\leq\sup_{s\in[t,T]}\abv{\Ebh\edg{\Yb^1_t}-\Ebh\edg{\Yb^2_t}}+\frac{2\cb}{\cu} \sup_{s\in[t,T]}\Ebh\edg{\abv{\Yb^1_s-\Yb_s^2}}\label{eq:DRapr}\\
		&\quad+\frac1\cu\Bigg\{\sup_{s\in[t,T]}\max\edg{\abv{\Hb^1(s,L^1_s(\Yb^2_s),\Yb_s^2)-\Hb^2(s,L^1_s(\Yb_s^2),\Yb_s^2)},\abv{\Hb^1(s,U^1_s(\Yb^2_s),\Yb_s^2)-\Hb^2(s,U^1_s(\Yb_s^2),\Yb_s^2)}}\notag\\
		&\quad\quad\quad\quad+\sup_{s\in[t,T]}\max(\abv{l^1_s-l_s^2},\abv{u^1_s-u^2_s})\Bigg\}\notag,
	\end{align}
    where the mappings $\Hb^i(s,\cdot,X),~L^i_s(X),~U^i_s(X),~X\in L_G^1(\Om_s),~s\in[0,T],~i=1,2$ are defined w.r.t. $\hbar^i,l^i,u^i$ correspondingly.
\end{Proposition}
\noindent{\bf\underline{Proof.}} For $i=1,2$ and $t\in[0,T]$, set 
\begin{align*}
	\tau_1^i:=0\vee\sup\crl{t<T:~\Ebh\edg{\hbar^i(t,Y^i_t)}=u^i_t~\mbox{or}~~\Ebh\edg{\hbar^i(t,Y^i_t)}=l_t^i},
\end{align*}
and for $k\geq2$,
\begin{align*}
	&\tau_k^i:=0\vee\sup\bigg\{t<\tau^i_{k-1}:~\Ebh\edg{\hbar^i(t,Y^i_t)}=u^i_t~\mbox{if}~~\Ebh\edg{\hbar^i(\tau^i_{k-1},Y^i_{\tau^i_{k-1}})}=l^i_{\tau^i_{k-1}},\\
		&~~~~~~~~~~~~~~~~~~~~~~~~~~~~~~~~~~\mbox{or}~~\Ebh\edg{\hbar^i(t,Y^i_t)}=l^i_t~\mbox{if}~~\Ebh\edg{\hbar^i(\tau^i_{k-1},Y^i_{\tau^i_{k-1}})}=u^i_{\tau^i_{k-1}}\bigg\},
\end{align*}
with the convention that $\sup\emptyset=-\infty$.
Denote $\Tc=\crl{\tau^i_k,i=1,2,k\in\N^*}$, $t_1:=0\vee\sup\crl{t:~t\in\Tc}$, and $t_{k+1}:=0\vee\sup\crl{t<t_k:~t\in\Tc},~k\in\N^*$. Similarly to Lemma 4.1 in \cite{Diffintreverse} and estimate \eqref{eq:infUL},  $$\eps^i=\inf_{t\in[0,T]}(U^i_t(\Yb^i_t)-L^i_t(\Yb^i_t))>0,~i=1,2.$$ 
There exist some $\delta^i>0$ such that,
\begin{align*}
	\abv{\Ebh\edg{Y^i_t}-\Ebh\edg{Y^i_s}}\leq\frac{\eps^i}{3},~\abv{L_t^i(\Yb^i_t)-L_s^i(\Yb^i_s)}\leq\frac{\eps^i}3,~\mbox{and}~~\abv{U^i_t(\Yb^i_t)-U^i_s(\Yb^i_s)}\leq\frac{\eps^i}{3}, 
\end{align*} 
whenever $\abv{s-t}\leq \delta^i$.
As a result, if $\tau_k^i>0$, $\tau^i_{k-1}-\tau_k^i>\delta^i$, which implies that $\tau^i_{N^i}=0$ for some positive integers $ N^i<+\infty$. Thus $\Tc$ is finite if not empty.

Next, we prove by induction inspired by Proposition 2.3 in \cite{SPtimebs}. For $t\in[t_1,T]$, the result is trivial. Assume for all $t\in[t_k,T]$, equality \eqref{eq:DRapr} holds. For $t\in[t_{k+1},t_k]$, we can split the problem into three cases:
\begin{enumerate}[(a)]
	\item exactly one of $\Ebh\edg{\hbar^1(\cdot,Y^1_\cdot)}$ and $\Ebh\edg{\hbar^2(\cdot,Y^2_\cdot)}$ hits a barrier at time $t_k$;
	\item both $\Ebh\edg{\hbar^1(\cdot,Y^1_\cdot)}$ and $\Ebh\edg{\hbar^2(\cdot,Y^2_\cdot)}$ hit upper barriers, or both $\Ebh\edg{\hbar^1(\cdot,Y^1_\cdot)}$ and $\Ebh\edg{\hbar^2(\cdot,Y^2_\cdot)}$ hit lower barriers at time $t_k$;
	\item $\Ebh\edg{\hbar^1(\cdot,Y^1_\cdot)}$ hits the lower barrier and $\Ebh\edg{\hbar^2(\cdot,Y^2_\cdot)}$ hits the upper barrier at time $t_{k}$, or $\Ebh\edg{\hbar^1(\cdot,Y^1_\cdot)}$ hits the upper barrier and $\Ebh\edg{\hbar^2(\cdot,Y^2_\cdot)}$ hits the lower barrier at time $t_{k}$.
\end{enumerate}
In view of Definition \ref{definition:BSPSE} and Remark \ref{remark: Rflat}, we concentrate on the most complicated case (c).  Assume w.l.o.g. $\Ebh\edg{Y^1_\cdot}$ hits $L^1_\cdot(\Yb^1_\cdot)$ and $\Ebh\edg{Y^2_\cdot}$ hits $U^2_\cdot(\Yb^2_\cdot)$ at time $t_k$.
Set 
\begin{equation*}
	\Yb^{i,*}_t:=\left\{
	\begin{lgathered}
		\Yb^i_t+R^i_T-R^i_{t_k},~t\in[t_{k+1},t_k];\\
		Y^i_t,~t\in(t_k,T],
	\end{lgathered}
	\right.
\end{equation*}
and $(Y^{i,*},R^{i,*}):=\widehat{\B\Sb\P}_{l^i}^{u^i}(\Yb^{i,*},\hbar^i)$. Observe that $$(Y^i_t,R^i_T-R^i_t)=(Y^{i,*}_t,R^{i,*}_{t_k}-R^{i,*}_t)+(0,R^i_T-R^i_{t_k}),~t\in[t_{k+1},t_k], ~i=1,2.$$
Applying Theorem \ref{theorem:BSPE&U}, we have
\begin{equation*}
	Y^{1,*}_t=\left\{
	\begin{lgathered}
		\Yb^1_t+R^1_T-R^1_{t_k}+\sup_{s\in[t,t_k]}\brak{\Ebh\edg{\Yb^1_s}+R^1_T-R^1_{t_k}-L_s^1(\Yb^1_s)}^-,~t\in[t_{k+1},t_k];\\
		Y^1_t,~t\in(t_k,T].
	\end{lgathered}
	\right.
\end{equation*}
For $t\in[t_{k+1},t_k]$, we have
\begin{equation*}
	\Yb^1_t+R^1_T-R^1_{t_k}+\sup_{s\in[t,t_k]}\brak{\Ebh\edg{\Yb^1_s}+R^1_T-R^1_{t_k}-L_s^1(\Yb^1_s)}^-=\Yb^1_t+\sup_{s\in[t,t_k]}\brak{L_s^1(\Yb^1_s)-\Ebh\edg{\Yb^1_s}}.
\end{equation*}
Similarly, we can get 
\begin{equation*}
	Y^{2,*}_t=\left\{
	\begin{lgathered}
			\Yb^2_t-\sup_{s\in[t,t_k]}\brak{\Ebh\edg{\Yb^2_s}-U^2_s(\Yb^2_s)},~t\in[t_{k+1},t_k];\\
				Y^2_t,~t\in(t_k,T].
	\end{lgathered}
	\right.
\end{equation*}
Clearly, we can obtain that for $t\in[t_{k+1},t_k]$,
$$(R^1_T-R^1_t)-(R^2_T-R^2_t)=\sup_{s\in[t,t_k]}\brak{L_s^1(\Yb^1_s)-\Ebh\edg{\Yb^1_s}}+\sup_{s\in[t,t_k]}\brak{\Ebh\edg{\Yb^2_s}-U^2_s(\Yb^2_s)}=:I_1(t)+I_2(t).$$
 Since $L^1_\cdot(\Yb^1_\cdot),\Ebh\edg{\Yb^1_\cdot},U^2_\cdot(\Yb^2_\cdot),\Ebh\edg{\Yb^2_\cdot}$ are continuous in $t$, there exist $t_1',t_2'\in[t,t_k]$ such that $I_1(t)=L^1_{t_1'}(\Yb^1_{t'_1})-\Ebh\edg{\Yb^1_{t'_1}}$ and $I_2(t)=\Ebh\edg{\Yb^2_{t'_2}}-U^2_{t_2'}(\Yb^2_{t_2'})$. 
 
 Suppose that $I_1(t)+I_2(t)>0$. If $t_1'\leq t'_2$, then
\begin{align*}
	 &\quad I_1(t)+I_2(t)\\
	 &=\Ebh\edg{\Yb^2_{t'_2}}-\Ebh\edg{\Yb^2_{t'_1}}+\Ebh\edg{\Yb^2_{t'_1}}-\Ebh\edg{\Yb^1_{t'_1}}+L^1_{t'_1}(\Yb^1_{t'_1})-U^2_{t'_2}(\Yb^2_{t'_2})\\
	&=\Ebh\edg{Y^2_{t'_2}}-(R^2_{T}-R^2_{t'_2})-\Ebh\edg{Y^2_{t'_1}}+(R^2_{T}-R^2_{t'_1})+\Ebh\edg{\Yb^2_{t'_1}}-\Ebh\edg{\Yb^1_{t'_1}}+L^1_{t'_1}(\Yb^1_{t'_1})-U^2_{t'_2}(\Yb^2_{t'_2})\\
	&\leq\Ebh\edg{\Yb^2_{t'_1}}-\Ebh\edg{\Yb^1_{t'_1}}+L^1_{t'_1}(\Yb^1_{t'_1})-L^2_{t'_1}(\Yb^2_{t'_1}),
\end{align*}
where the last inequality is deduced by the fact that $dR^2_s=0,~s\in[t'_1,t'_2]$, $\Ebh\edg{Y^2_{t'_2}}\leq U^2_{t'_2}(\Yb^2_{t'_2})$, and $\Ebh\edg{Y^2_{t'_1}}\geq L_{t'_1}^2(\Yb^2_{t'_1})$.
If $t'_2<t'_1$, then,
\begin{align*}
 &\quad I_1(t)+I_2(t)\\
 &=\Ebh\edg{\Yb^1_{t'_2}}-\Ebh\edg{\Yb^1_{t'_1}}+\Ebh\edg{\Yb^2_{t'_2}}-\Ebh\edg{\Yb^1_{t'_2}}+L^1_{t'_1}(\Yb^1_{t'_1})-U^2_{t'_2}(\Yb^2_{t'_2})\\
	&=\Ebh\edg{Y^1_{t'_2}}-(R^1_{T}-R^1_{t'_2})-\Ebh\edg{Y^1_{t'_1}}+(R^1_{T}-R^1_{t'_1})+\Ebh\edg{\Yb^2_{t'_2}}-\Ebh\edg{\Yb^1_{t'_2}}+L^1_{t'_1}(\Yb^1_{t'_1})-U^2_{t'_2}(\Yb^2_{t'_2})\\
	&\leq\Ebh\edg{\Yb^2_{t'_2}}-\Ebh\edg{\Yb^1_{t'_2}}+U^1_{t'_2}(\Yb^1_{t'_2})-U^2_{t'_2}(\Yb^2_{t'_2}).
\end{align*}
Note that for any $s\in[0,T]$, it holds that
\begin{align*}
	&\quad\abv{L^1_s(\Yb^1_s)-L^2_s(\Yb^2_s)}\leq \abv{L^1_s(\Yb^1_s)-L^1_s(\Yb^2_s)}+\abv{L^1_s(\Yb^2_s)-L^2_s(\Yb^2_s)}\\
	&\leq\frac{2\cb}{\cu}\Ebh\edg{\abv{\Yb^1_s-\Yb_s^2}}+\frac1\cu\abv{\Ebh\edg{\hbar^2(s,L^1_s(\Yb^2_s)+\Yb^2_s-\Ebh\edg{\Yb^2_s})}-\Ebh\edg{\hbar^2(s,L^2_s(\Yb^2_s)+\Yb^2_s-\Ebh\edg{\Yb^2_s})}}\\
	&\leq\frac{2\cb}{\cu}\Ebh\edg{\abv{\Yb^1_s-\Yb_s^2}}+\frac1\cu\abv{\Ebh\edg{\hbar^1(s,L^1_s(\Yb^2_s)+\Yb^2_s-\Ebh\edg{\Yb^2_s})}-\Ebh\edg{\hbar^2(s,L^1_s(\Yb^2_s)+\Yb^2_s-\Ebh\edg{\Yb^2_s})}}\\
	&\quad+\frac1\cu\abv{\Ebh\edg{\hbar^1(s,L^1_s(\Yb^2_s)+\Yb^2_s-\Ebh\edg{\Yb^2_s})}-\Ebh\edg{\hbar^2(s,L^2_s(\Yb^2_s)+\Yb^2_s-\Ebh\edg{\Yb^2_s})}}\\
	&\leq\frac{2\cb}{\cu}\Ebh\edg{\abv{\Yb^1_s-\Yb_s^2}}+\frac1\cu\brak{\abv{\Hb^1(s,L^1_s(\Yb^2_s),\Yb^2_s)-\Hb^2(s,L^1_s(\Yb^2_s),\Yb^2_s)}+\abv{l^1_s-l^2_s}},
\end{align*}
and,
\begin{align*}
	\abv{U^1_s(\Yb^1_s)-U^2_s(\Yb^2_s)}\leq\frac{2\cb}{\cu}\Ebh\edg{\abv{\Yb^1_s-\Yb_s^2}}+\frac1\cu\brak{\abv{\Hb^1(s,U^1_s(\Yb^2_s),\Yb^2_s)-\Hb^2(s,U^1_s(\Yb^2_s),\Yb^2_s)}+\abv{u^1_s-u^2_s}}.
\end{align*}

Suppose that $I_1(t)+I_2(t)\leq0$. As $I_1(t)+I_2(t)$ is non-increasing in $t$, we can deduce that
\begin{align*}
	&L^1_{t_k}(\Yb^1_{t_k})-\Ebh\edg{\Yb^1_{t_k}}+\Ebh\edg{\Yb^2_{t_k}}-U^2_{t_k}(\Yb^2_{t_k})=(R^1_T-R^1_{t_k})-(R^2_T-R^2_{t_k})\leq I_1(t)+I_2(t)\leq0.
\end{align*}
Hence,
$
	~\abv{I_1(t)+I_2(t)}=\abv{(R^1_T-R^1_t)-(R^2_T-R^2_t)}\leq\abv{(R^1_T-R^1_{t_k})-(R^2_T-R^2_{t_k})}.
$

To summarize, for any $t\in[t_{k+1},t_k]$,
\begin{align*}
	&\quad\abv{(R^1_T-R^1_t)-(R^2_T-R^2_t)}\\
	&\leq \max\Bigg\{\sup_{s\in[t,t_k]}\abv{\Ebh\edg{\Yb^1_s}-\Ebh\edg{\Yb^2_s}}+\sup_{s\in[t,t_k]}\max\brak{\abv{L^1_s(\Yb_s^1)-L^2_s(\Yb^2_s)},\abv{U^1_s(\Yb^1_s)-U^2_s(\Yb^2_s)}},\\
	&~~~~~~~~~~~~~\abv{(R^1_T-R^1_{t_k})-(R^2_T-R^2_{t_k})}\Bigg\}\\
	&\leq \max\Bigg\{\sup_{s\in[t,t_k]}\abv{\Ebh\edg{\Yb^1_s}-\Ebh\edg{\Yb^2_s}}+\frac{2\cb}{\cu}\sup_{s\in[t,t_k]}\Ebh\edg{\abv{\Yb^1_s-\Yb_s^2}}\\
	&~~~~~~~~~~~~~~~~~+\frac1\cu\bigg\{\sup_{s\in[t,t_k]}\max\bigg[\abv{\Hb^1(s,L^1_s(\Yb^2_s),\Yb^2_s)-\Hb^2(s,L^1_s(\Yb^2_s),\Yb^2_s)},\\
	&~~~~~~~~~~~~~~~~~~~~~~~~~~\abv{\Hb^1(s,U^1_s(\Yb^2_s),\Yb^2_s)-\Hb^2(s,U^1_s(\Yb^2_s),\Yb^2_s)}\bigg]+\sup_{s\in[t,t_k]}\max\brak{\abv{l^1_s-l^2_s},\abv{u^1_s-u^2_s}}\bigg\},\\
	&~~~~~~~~~~~~~\sup_{s\in[t_k,T]}\abv{\Ebh\edg{\Yb^1_s}-\Ebh\edg{\Yb^2_s}}+\frac{2\cb}{\cu}\sup_{s\in[t_k,T]}\Ebh\edg{\abv{\Yb^1_s-\Yb_s^2}}\\
	&~~~~~~~~~~~~~~~~~+\frac1\cu\bigg\{\sup_{s\in[t_k,T]}\max\bigg[\abv{\Hb^1(s,L^1_s(\Yb^2_s),\Yb^2_s)-\Hb^2(s,L^1_s(\Yb^2_s),\Yb^2_s)},\\
	&~~~~~~~~~~~~~~~~~~~~~~~~~~\abv{\Hb^1(s,U^1_s(\Yb^2_s),\Yb^2_s)-\Hb^2(s,U^1_s(\Yb^2_s),\Yb^2_s)}\bigg]+\sup_{s\in[t_k,T]}\max\brak{\abv{l^1_s-l^2_s},\abv{u^1_s-u^2_s}}\bigg\}\Bigg\}\\
	&\leq\sup_{s\in[t,T]}\abv{\Ebh\edg{\Yb^1_s}-\Ebh\edg{\Yb^2_s}}+\frac{2\cb}{\cu}\sup_{s\in[t,T]}\Ebh\edg{\abv{\Yb^1_s-\Yb_s^2}}\\
	&~~~~~~~~~~~~~~~~~+\frac1\cu\Bigg\{\sup_{s\in[t,T]}\max\bigg[\abv{\Hb^1(s,L^1_s(\Yb^2_s),\Yb^2_s)-\Hb^2(s,L^1_s(\Yb^2_s),\Yb^2_s)},\\
	&~~~~~~~~~~~~~~~~~~~~~~~~~~\abv{\Hb^1(s,U^1_s(\Yb^2_s),\Yb^2_s)-\Hb^2(s,U^1_s(\Yb^2_s),\Yb^2_s)}\bigg]+\sup_{s\in[t,T]}\max\brak{\abv{l^1_s-l^2_s},\abv{u^1_s-u^2_s}}\Bigg\}.
\end{align*}
The proof is therefore ended.
\qed

\begin{Proposition}\label{proposition:BSPSEapr}
	Let Assumption \ref{assumption:HR} hold, and $\Yb\in C_G^1(0,T)$ with $l_T\leq\Ebh\edg{\hbar(T,\Yb_T)}\leq u_T$. If $C_G^1(0,T)\x C_B^T\ni(Y,R)=\widehat{\B\Sb\P}_{l}^{u}(\Yb,\hbar)$, then
	\begin{align}
		\abv{R_T-R_t}\leq \Cs \sup_{s\in[t,T]}\Ebh\edg{\abv{\Yb_s}}+C_0(t,T),\label{eq:BSPYapr}
	\end{align}
   where $\Cs:=1+\frac{2\cb}\cu$, and $C_0(t,T):=\sup_{s\in[t,T]}\max\edg{\abv{L_s(0)},\abv{U_s(0)}}$.
\end{Proposition}
\noindent{\bf\underline{Proof.}} 
Similarly to Proposition \ref{proposition:BSPSEaprDf}, we set
\begin{align*}
	\tau_1:=0\vee\sup\crl{t<T:~\Ebh\edg{\hbar(t,Y_t)}=u_t ~\mbox{or}~~\Ebh\edg{\hbar(t,Y_t)}=l_t},
\end{align*}
and for $k\geq2$,
\begin{align*}
	&\tau_k:=0\vee\sup\Bigg\{t<\tau_k:~\Ebh\edg{\hbar(t,Y_t)}=u_t~\mbox{if}~~\Ebh\edg{\hbar(\tau_{k-1},Y_{\tau_{k-1}})}=l_{\tau_{k-1}},\\
	&~~~~~~~~~~~~~~~~~~~~~~~~~~~~~~~\mbox{or}~~\Ebh\edg{\hbar(t,Y_t)}=l_t~\mbox{if}~~\Ebh\edg{\hbar(\tau_{k-1},Y_{\tau_{k-1}})}=u_{\tau_{k-1}}	\Bigg\}.
\end{align*}
Accordingly, there exists some $\N^*\ni N_0<+\infty$ such that $\tau_{N_0}=0$. 
Note that on $[\tau_1,T]$,  inequality \eqref{eq:BSPYapr} holds trivially.
Suppose that inequality
\eqref{eq:BSPYapr} holds on $t\in[\tau_k,T]$.
Assume that $\Ebh\edg{\hbar(\cdot,Y_{\cdot})}$ hits $l$ at time $\tau_k$. Put
\begin{equation*}
	\Yb^*_t=\left\{
	\begin{lgathered}
		\Yb_t+R_T-R_{\tau_k},~t\in[\tau_{k+1}, \tau_k];\\
		Y_t,~t\in(\tau_k,T],
	\end{lgathered}
	\right.
\end{equation*} 
and $(Y^*,R^*)=\widehat{\B\Sb\P}_l^u(\Yb^*,\hbar)$. Hence $(Y_t,R_T-R_t)=(Y^*_t,R^*_{\tau_k}-R^*_t)+(0,R_T-R_{\tau_k}),~t\in[\tau_{k+1}, \tau_k]$. By Theorem \ref{theorem:BSPE&U}, we have on $t\in[\tau_{k+1},\tau_k]$,
\begin{equation*}
	Y^*_t=\Yb_t+R_T-R_{\tau_k}+\sup_{s\in[t,\tau_k]}\brak{\Ebh\edg{\Yb_s}+R_T-R_{\tau_k}-L_s(\Yb_s)}^-=\Yb_t+\sup_{s\in[t,\tau_k]}\brak{L_s(\Yb_s)-\Ebh\edg{\Yb_s}}.
\end{equation*}
If $R_T-R_t=\sup_{s\in[t,\tau_k]}\brak{L_s(\Yb_s)-\Ebh\edg{\Yb_s}}\geq0$, then
\begin{align*}
	\abv{R_T-R_t}\leq \sup_{s\in[t,\tau_k]}\Ebh\edg{\abv{\Yb_s}}+\sup_{s\in[t,\tau_k]}\abv{L_s(\Yb_s)}\leq (1+\frac{2\cb}{\cu})\sup_{s\in[t,\tau_k]}\Ebh\edg{\abv{\Yb_s}}+\sup_{s\in[t,\tau_k]}\abv{L_s(0)}.
\end{align*}
If $R_T-R_t<0$, then
$
	R_T-R_{\tau_k}\leq R_T-R_t<0,~\mbox{and}~~\abv{R_T-R_t}\leq \abv{R_T-R_{\tau_k}}$.
Therefore,
\begin{align*}
	\abv{R_T-R_t}
	\leq (1+\frac{2\cb}{\cu}) \sup_{s\in[t,T]}\Ebh\edg{\abv{\Yb_s}}+\sup_{s\in[t,T]}\max\edg{L_s(0),U_s(0)}.
\end{align*}
We end the proof by mentioning that the results where $\Ebh\edg{\hbar(\cdot,Y_\cdot)}$ hits $u$ at time $\tau_k$ can be achieved analogously.
\qed

\section{Doubly mean-reflected $G$-BSDEs}
In this section, we consider a series of DMR$G$-BSDEs $(\xi,f,\hbar,l,u)$:
\begin{equation}
	\left\{
	\begin{lgathered}
		Y_t=\xi+\int_t^Tf(s,Y_s,Z_s)ds-\int_t^TZ_sdB_s-(K_T-K_t)+(R_T-R_t),~t\in[0,T];\\
		l_t\leq\Ebh\edg{\hbar(t,Y_t)}\leq u_t,~\forall t\in[0,T];\label{eq:mainDMR}\\
	\int_{s}^{t}\brak{\Ebh\edg{\hbar(v,Y_v)}-l_v}dR_v\leq0,~\mbox{and}~\int_{s}^{t}\brak{\Ebh\edg{\hbar(v,Y_v)}-u_v}dR_v\leq0,~ \forall0\leq s\leq t\leq T;
	\end{lgathered}
	\right.
\end{equation} 
where a universal hypothesis holds for the terminal condition:
\begin{Assumption}\label{assumption:xireflection} $\xi\in L_G^1(\Om_T)$ and
	$l_T\leq \Ebh\edg{\hbar(T,\xi)}\leq u_T$.
\end{Assumption}
To establish corresponding well-posedness results, using backward Skorokhod problems with sublinear expectation, we perform briefly fixed-point approaches and $\theta$-methods inspired by works on standard $G$-BSDEs, see \cite{GLX, QGBSDEbt,QGBSDEubt,LiumdGBSDE,LWMFGB} for more details.
\subsection{Lipschitz generators}

\begin{Assumption}\label{assumption:HLip}
Fix some constants $\b>1$ and $L>0$.
	\begin{enumerate}[($i$)]
		\item $\xi\in L_G^\b(\Om_T)$.
		\item For any $(t,\om)\in[0,T]\x\Om_T$, $f(t,\om,0,\0)\in M_G^\b(0,T)$.
		\item For each $(t,\om)\in[0,T]\x\Om_T$ and any $(y,y',z,z')\in\R^{2+2d}$,
		\begin{equation*}
			\abv{f(t,\om,y,z)-f(t,\om,y',z')}\leq L\brak{\abv{y-y'}+\abv{z-z'}}.
		\end{equation*}
	\end{enumerate}
\end{Assumption}

\begin{Theorem}\label{theorem:LipDMR}
	Let Assumptions \ref{assumption:HR}-\ref{assumption:xireflection}-\ref{assumption:HLip} hold. Then for any $1\leq \a<\b$, DMR$G$-BSDE $(\xi,f,\hbar,l,u)$ \eqref{eq:mainDMR} admits a unique solution $(Y,Z,K,R)\in\Sc_G^\a(0,T)\x C_B^T$.
\end{Theorem}
\noindent{\bf\underline{Proof.}} Take $U^{i}\in M_G^\b(0,T)$, $i=1,2$. Then the following $G$-BSDE:
\begin{equation}\label{eq:testGBSDE}
	\Yb^{i}=\xi+\int_t^Tf(s,U^{i}_s,\Zb^{i})ds-\int_t^T\Zb^{i}dB_s-(\Kb^{i}_T-\Kb^{i}_t),~t\in[0,T],
\end{equation}
admits a unique solution $(\Yb^{L,i},\Zb^{L,i},\Kb^{L,i})\in\Sc_G^\a(0,T)$, for any $1\leq\a<\b$, from Theorem 4.1 in \cite{HuJiPengSong}. By Theorem \ref{theorem:BSPE&U}, there exist unique solutions to backward Skorokhod problem with sublinear expectation, i.e., $C_G^1(0,T)\x C_B^T\ni (Y^{L,i},R^{L,i})=\widehat{\B\Sb\P}_l^u(\Yb^{L,i},\hbar))$, $i=1,2$. Therefore, putting $f^{L,i}(s,z):=f(s,U^{i}_s,z)$,
$(Y^{L,i},\Zb^{L,i},\Kb^{L,i},R^{L,i})$
satisfies DMR$G$-BSDE \eqref{eq:mainDMR} ($\xi,f^{L,i},\hbar,l,u$), uniquely.  By Proposition \ref{proposition:BSPSEapr}, $Y^{L,i}$ share the same integrability with $\Yb^{L,i}$, and they all belong to $M_G^\b(0,T)$ by Lemma 3.2 in \cite{LiumdGBSDE}.
According to Proposition \ref{prop:1daprDY}, there exists a constant $C^L:=C^L(T,\sigmat,\sigmah,L,\b)>0$ such that for $t\in[0,T]$,
\begin{align*}
	&\abv{\Yb^{L,1}_t-\Yb^{L,2}_t}^\b\\
	&\leq C^L\Ebh_t\edg{\brak{\int_t^T\abv{f(s,U^{1}_s,\Zb^{L,2})-f(s,U^{2}_s,\Zb^{L,2}_s)}ds}^\b}\leq C^LL^\b(T-t)^{\b-1}\Ebh_t\edg{\int_t^T\abv{U^{1}_s-U^{2}_s}^\b ds},
\end{align*} 
and
\begin{align*}
	\Ebh\edg{\int_t^T\abv{\Yb^{L,1}_s-\Yb^{L,2}_s}^\b ds}\leq \int_t^T\Ebh\edg{\abv{\Yb^{L,1}_s-\Yb^{L,2}_s}^\b}ds\leq C^LL^\b(T-t)^\b\Ebh\edg{\int_t^T\abv{U^{1}_s-U^{2}_s}^\b ds}.
\end{align*}
Concerning Proposition \ref{proposition:BSPSEaprDf} and $\Cs=1+\frac{2M}{m}$, we have
\begin{align*}
\Ebh\edg{\int_t^T\abv{Y^{L,1}_s-Y^{L,2}_s}^\b ds}
&\leq \int_t^T\Ebh\edg{\abv{Y^{L,1}_s-Y^{L,2}_s}^\b}ds\\
&\leq \int_t^T\Ebh\edg{\brak{\abv{\Yb^{L,1}_s-\Yb^{L,2}_s}+\Cs\sup_{v\in[s,T]}\Ebh\edg{\abv{\Yb^{L,1}_v-\Yb^{L,2}_v}}}^\b}ds\\
&\leq 2^{\b-1}\int_t^T\Ebh\edg{\abv{\Yb^{L,1}_s-\Yb^{L,2}_s}^\b}ds+2^{\b-1}\Cs^\b(T-t)\sup_{s\in[t,T]}\Ebh\edg{\abv{\Yb^{L,1}_s-\Yb^{L,2}_s}^\b}ds\\
&\leq 2^{\b-1}(1+\Cs^\b)C^LL^\b(T-t)^\b\Ebh\edg{\int_t^T\abv{U^{1}_s-U^{2}_s}^\b ds}.
\end{align*}
Taking $h_1\in(0,T]$ such that $2^{\b-1}(1+\Cs^\b)C^LL^\b h_1^\b<1$ and $T=N_h^1h_1$ for some $N_h^1\in\N^*$, we can construct a contraction mapping on $M_G^\b(T-h_1,T)$ to derive a fixed point $Y\in M_G^\b(T-h_1,T)$ and then local solutions by means of a priori estimates for $Z$, see Proposition \ref{prop:1daprDZ}. Iterating backwardly, we get the desired convergence and uniqueness results. 
\qed

\subsection{Quadratic generators}
\subsubsection{Bounded terminal conditions}

\begin{Assumption}\label{assumption:HQbdd}Let $\lm,\gamma,M_0\in\R^+$.
	\begin{enumerate}[($i$)]
		\item $\xi\in L_G^\infty(\Om_T)$.
		\item There exists some constant $M_0>0$ such that for each $(t,\om)\in[0,T]\x\Om_T$,
		$$\abv{\xi}+\int_0^T\abv{f(s,\om,0,\0)}^2ds\leq M_0.$$
		\item For any $(t,\om)\in[0,T]\x\Om_T$ and any  $(y,y',z,z')\in\R^{2+2d}$, 
		$$\abv{f(t,\om,y,z)-f(t,\om,y',z')}\leq\lm\abv{y-y'}+\gamma(1+\abv{z}+\abv{z'})\abv{z-z'}.$$
		\item There is a modulus of continuity $w:[0,\infty)\to[0,\infty)$ such that
		$$\abv{f(t,\om,y,z)-f(t',\om',y,z)}\leq w(\abv{t-t'}+\norm{\om-\om'}),$$
		for each $(y,z)\in\R^{1+d}$ and any  $(t,t',\om,\om')\in[0,T]\x[0,T]\x\Om_T\x\Om_T$.
	\end{enumerate}
\end{Assumption}

\begin{Theorem}\label{theorem:btv-2}
	Let Assumptions \ref{assumption:HR}-\ref{assumption:xireflection}-\ref{assumption:HQbdd} hold. Then DMR$G$-BSDE $(\xi,f,\hbar,l,u)$ \eqref{eq:mainDMR} 
	admits a unique solution $(Y,Z,K,R)\in S_G^\infty(0,T)\x \Hs_G^T(\R^d)\x\As_G^T\x C_B^T$.
\end{Theorem}
\noindent{\bf\underline{Proof.}} Take $U^{i}\in S_G^\infty(0,T),~i=1,2$ and then the Q$G$-BSDE \eqref{eq:testGBSDE}
admit a unique solution $(\Yb^{Q,i},\Zb^{Q,i},\Kb^{Q,i})\in S_G^\infty(0,T)\x\Hs_G^T(\R^d)\x\As_G^T$ by Remark \ref{remark:QGZ}.  Then applying Theorem \ref{theorem:BSPE&U}, there exists uniquely $C_G^1(0,T)\x C_B^T\ni(Y^{Q,i},R^{Q,i})=\widehat{\B\Sb\P}_l^u(\Yb^{Q,i},\hbar)$.
Therefore, putting $f^{Q,i}(s,z):=f(s,U_s^{Q,i},z)$, $(Y^{Q,i},\Zb^{Q,i},\Kb^{Q,i},R^{Q,i})$ satisfies DMR$G$-BSDE ($\xi,f^{Q,i},\hbar,l,u$) \eqref{eq:mainDMR} uniquely. And from Proposition \ref{proposition:BSPSEapr}, $Y^{Q,i}\in S_G^\infty(0,T)$.

Following a similar $G$-Girsanov transform approach to Theorem 3.2 of \cite{GLX}, we derive that,
\begin{align*}
	\norm{\Yb^{Q,1}-\Yb^{Q,2}}_{S_G^\infty(t,T)}\leq \lm (T-t)\norm{U^{Q,1}-U^{Q,2}}_{S_G^\infty(t,T)},~t\in[0,T]. 
\end{align*}
Referring to Proposition \ref{proposition:BSPSEaprDf}, we deduce that
\begin{align*}
	\norm{Y^{Q,1}-Y^{Q,2}}_{S_G^\infty(t,T)}\leq(1+\Cs)\norm{\Yb^{Q,1}-\Yb^{Q,2}}_{S_G^\infty(t,T)}\leq  (1+\Cs)\lm (T-t)\norm{U^{Q,1}-U^{Q,2}}_{S_G^\infty(t,T)}.
\end{align*}
Choose $h_2\in(0,T]$ such that $(1+\Cs)\lm h_2<1$ and $T=N^{2}_hh_2$ for some $N^{2}_h\in\N^*$. Applying Proposition \ref{proposition:qgbsdeaprdyz}, we can employ a fixed-point argumentation to derive local solutions and concatenate them afterwards. And the uniqueness is also provided by local results.
\qed 

\subsubsection{Unbounded terminal conditions}

\begin{Assumption}\label{assumption:HQUbdd}
Let $\lm,\gamma\in\R^+$ and $\ha\in M_G^1(0,T)$ be nonnegative.
	\begin{enumerate}[($i$)]
		\item Both $\xi$ and $\int_0^T\ha_tdt$ have finite exponential moments of arbitrary order, i.e.,
		$$\Ebh\edg{\exp\crl{p\abv{\xi}+p\int_0^T\ha_tdt}}<\infty,~\forall p\geq1.$$
		\item For any $(t,\om)\in[0,T]\x\Om_T$ and any  $(y,y',z,z')\in\R^{2+2d}$, 
		$$\abv{f(t,\om,y,z)-f(t,\om,y',z')}\leq\lm\abv{y-y'}+\gamma(1+\abv{z}+\abv{z'})\abv{z-z'}.$$
		\item There is a modulus of continuity $w:[0,\infty)\to[0,\infty)$ such that
		$$\abv{f(t,\om,y,z)-f(t',\om',y,z)}\leq w(\abv{t-t'}+\norm{\om-\om'}),$$
		for each $(y,z)\in\R^{1+d}$ and any  $(t,t',\om,\om')\in[0,T]\x[0,T]\x\Om_T\x\Om_T$.
		\item For each $(t,\om,y)\in[0,T]\x\Om_T\x\R$, $f(t,\om,y,\cdot)$ is convex.
\end{enumerate}
\end{Assumption}

Note that for any $(t,\om,y,z)\in[0,T]\x\Om_T\x\R\x\R^d$,
\begin{equation}\label{eq:fbound}
	\abv{f(t,\om,y,z)}\leq\abv{f(t,\om,0,\0)}+\lm\abv{y}+\gamma(\abv{z}^2+\abv{z})\leq\ha_t(\om)+\frac{\gamma}{2}+\lm\abv{y}+\frac{3\gamma}{2}\abv{z}^2.
\end{equation} 	

\begin{Theorem}\label{theorem:main}
	Let Assumptions \ref{assumption:HR}-\ref{assumption:xireflection}-\ref{assumption:HQUbdd} hold. Then DMRQ$G$-BSDE $(\xi,f,\hbar,l,u)$ \eqref{eq:mainDMR} admits a unique solution $(Y,Z,K,R)\in\Es_G^T\x\Hs_G^T(\R^d)\x\As_G^T\x C_B^T$.
\end{Theorem}

In order to apply a $\theta$-method following \cite{GLX,QGBSDEubt}, we denote for each $\theta\in(0,1)$,
\begin{equation*}
	\dt\phi:=\frac{\phi^1-\theta\phi^2}{1-\theta},~\dtt\phi:=\frac{\phi^2-\theta\phi^1}{1-\theta}~\mbox{and}~\dta\phi:=\abv{\dt\phi}+\abv{\dtt\phi},~\mbox{for}~~\phi=\Yb,Y,\Zb,R_T-R_\cdot,
\end{equation*}
and the following lemma is beneficial for proving required boundedness results. 
\begin{Lemma}\label{lemma:Rtheta}
	Let Assumptions \ref{assumption:HR} hold, $\theta\in(0,1)$, and $\Yb^i\in C_G^1(0,T)$ with $l_T\leq\Ebh\edg{\hbar(T,\Yb^i_T)}\leq u_T]$. If $C_G^1(t,T)\x C_B(t,T)\ni(Y^i,R^i)=\widehat{\B\Sb\P}_l^u(\Yb^i,\hbar)$, $i=1,2$. Then for $t\in[0,T]$,
	\begin{equation*}
		\abv{\dt Y_t}\leq C_0'(t,T,\Yb^1,\Yb^2)+\abv{\dt\Yb_t}+\Cs\sup_{s\in[t,T]}\Ebh\edg{\abv{\dt\Yb_s}},
	\end{equation*}
	and
	\begin{equation*}
		\abv{\dtt Y_t}\leq C_0'(t,T,\Yb^1,\Yb^2)+\abv{\dtt\Yb_t}+\Cs\sup_{s\in[t,T]}\Ebh\edg{\abv{\dtt\Yb_s}},
	\end{equation*}
	where
	\begin{align*}
		C_0'(t,T,\Yb^1,\Yb^2)   	&:=\frac2{\cu}\sup_{s\in[t,T]}\max\brak{\abv{l_s},\abv{u_s}}+\frac{\cb}{\cu}\sup_{s\in[t,T]}\max(\abv{L_s(0)},\abv{U_s(0)})\\
		&~~~~~~~~~~~~~~~~~~~~~~~~~~~~~~~~+\frac{2\cb(\cu+\cb)}{\cu^2}\brak{\sup_{s\in[t,T]}\Ebh\edg{\abv{\Yb^1_s}}+\sup_{s\in[t,T]}\Ebh\edg{\abv{\Yb^2_s}}},
	\end{align*}
	and as Proposition \ref{proposition:BSPSEapr},
	$
	\Cs:=1+\frac{2\cb}{\cu}.
	$
\end{Lemma}
\noindent{\bf\underline{Proof.}}
We reformulate the backward Skorokhod problems with sublinear expectation as $(\Yu^i,\Ru^i):=\widehat{\B\Sb\P}_{l^i}^{u^i}(\Ybu^i,\hbu^i),~i=1,2$ as: for each $t\in[0,T]$,
\begin{align*}
	\Yu^1_t:=Y^1_t,~\Ru^1_t:=R^1_t,~\hbu^1(t,\cdot):=\hbar(t,\cdot),~\Ybu^1_t:=\Yb^1_t,~l^1_t:=l_t,~u^1_t:=u_t;
\end{align*}
and 
\begin{align*}
	\Yu^2_t:=\theta Y^2_t,~\Ru^2_t:=\theta R^2_t,~\Ybu^2_t:=\theta\Yb^2_t, ~\hbu^2(t,y):=\theta\hbar(t,\frac1\theta y),~y\in\R,~l^2_t:=\theta l_t,~u^2_t:=\theta u_t.
\end{align*}
Further, putting  $\Hbu^i(t,x,\Ybu^i_t):=\Ebh\edg{\hbu^i(t,x+\Ybu^i_t-\Ebh\edg{\Ybu^i_t})},~i=1,2$, and defining correspondingly $\Lu^i_\cdot(\Ybu^i_\cdot),\Uu_\cdot^i(\Ybu^i_\cdot)$, Propositions \ref{proposition:hbarpropoty}-\ref{proposition:ULcon} hold.
By Proposition \ref{proposition:BSPSEaprDf}, we have,
\begin{align*}
	\abv{\dt Y_t}
	&\leq \abv{\dt\Yb_t}+\abv{\dt(R_T-R_t)}\leq  \abv{\dt\Yb_t}+\frac1{1-\theta}\abv{(\Ru^1_T-\Ru^2_t)-(\Ru^2_T-\Ru_t^2)}\\
	&\leq \abv{\dt\Yb_t}+ (1+\frac{2\cb}{\cu})\sup_{s\in[t,T]}\Ebh\edg{\abv{\dt\Yb_s}}+\frac1{\cu}\sup_{s\in[t,T]}\max\brak{\abv{l_s},\abv{u_s}}\\
	&\quad+\frac{1}{\cu(1-\theta)}\sup_{s\in[t,T]}\max\bigg[\abv{\Hbu^1(s,\Lu^1_s(\Ybu^2_s),\Ybu_s^2)-\Hbu^2(s,\Lu^1_s(\Ybu^2_s),\Ybu^2_s)},\\
	&~~~~~~~~~~~~~~~~~~~~~~~~~~~~~~~~~~~~~~~~~~~~~~~~~~~~~~~~~~\abv{\Hbu^1(s,\Uu^1_s(\Ybu^2_s),\Ybu_s^2)-\Hbu^2(s,\Uu^1_s(\Ybu^2_s),\Ybu^2_s)}\bigg].
\end{align*}
Here, we can derive separately that
\begin{align*}
	&\quad\abv{\Hbu^1(s,\Lu^1_s(\Ybu^2_s),\Ybu_s^2)-\Hbu^2(s,\Lu^1_s(\Ybu^2_s),\Ybu^2_s)}\\
	&=\abv{\Ebh\edg{\hbar(s,L_s(\theta\Yb^2_s)+\theta\Yb^2_s-\theta\Ebh\edg{\Yb^2_s})}-\Ebh\edg{\theta\hbar(s,\frac1\theta L_s(\theta\Yb^2_s)+\Yb^2_s-\Ebh\edg{\Yb^2_s})}}\\
	&\leq (1-\theta)\abv{l_s}+\theta\Ebh\edg{\abv{\hbar(s,L_s(\theta\Yb^2_s)+\theta\Yb^2_s-\theta\Ebh\edg{\Yb^2_s})-\hbar(s,\frac1\theta L_s(\theta\Yb^2_s)+\Yb^2_s-\Ebh\edg{\Yb^2_s})}}\\
	&\leq (1-\theta)\abv{l_s}+(1-\theta)\cb\abv{L_s(0)}+2\theta(1-\theta)\frac{\cb^2}{\cu}\Ebh\edg{\abv{\Yb^2_s}}+2\theta(1-\theta)\cb\Ebh\edg{\abv{\Yb^2_s}},
\end{align*}
and
\begin{align*}
	&\quad\abv{\Hbu^1(s,\Uu^1_s(\Ybu^2_s),\Ybu_s^2)-\Hbu^2(s,\Uu^1_s(\Ybu^2_s),\Ybu^2_s)}\\
	&\leq (1-\theta)\abv{u_s}+(1-\theta)\cb\abv{U_s(0)}+2\theta(1-\theta)\frac{\cb^2}{\cu}\Ebh\edg{\abv{\Yb^2_s}}+2\theta(1-\theta)\cb\Ebh\edg{\abv{\Yb^2_s}}.
\end{align*}
Consequently, we get
\begin{align*}
	\abv{\dt Y_t}
	&\leq \abv{\dt\Yb_t}+ (1+\frac{2\cb}{\cu})\sup_{s\in[t,T]}\Ebh\edg{\abv{\dt\Yb_s}}+\frac1{\cu}\sup_{s\in[t,T]}\max\brak{\abv{l_s},\abv{u_s}}\\
	&\quad+\frac{1}{\cu}\crl{\sup_{s\in[t,T]}\max(\abv{l_s},\abv{u_s})+\cb\sup_{s\in[t,T]}\max(\abv{L_s(0)},\abv{U_s(0)})+(\frac{2\cb^2}{\cu}+2\cb)\sup_{s\in[t,T]}\Ebh\edg{\abv{\Yb^2_s}}}\\
	&\leq \abv{\dt\Yb_t}+ (1+\frac{2\cb}{\cu})\sup_{s\in[t,T]}\Ebh\edg{\abv{\dt\Yb_s}}+\frac2{\cu}\sup_{s\in[t,T]}\max\brak{\abv{l_s},\abv{u_s}}+\frac{\cb}{\cu}\sup_{s\in[t,T]}\max(\abv{L_s(0)},\abv{U_s(0)})\\
	&\quad+\frac{2\cb(\cu+\cb)}{\cu^2}\sup_{s\in[t,T]}\Ebh\edg{\abv{\Yb^2_s}}.
\end{align*}
By symmetry, we can derive the estimate for $ 
\abv{\dtt Y_t}$.
\qed

\vspace{1em}
\noindent{\bf\underline{Proof of Theorem \ref{theorem:main}}}
  To start via a truncation argument, 
we denote by $\fb_t(\om):=f(t,\om,0,\0)$ and set for each $m\in\N^*$,
\begin{equation*}
	\phi^{(m)}:=(\phi\wedge m)\vee(-m)~\mbox{for}~\phi=\xi,\fb_t,~\mbox{and}~~f^{(m)}(t,y,z):=f(t,y,z)-\fb_t+\fb^{(m)}_t,
\end{equation*}
where $(\xi^{(m)},f^{(m)})$ satisfy Assumption \ref{assumption:HQbdd} with $M_0:=m+m^2T$, $\lm,\gamma>0$.
From Theorem \ref{theorem:btv-2}, DMR$G$-BSDE $(\xi^{(m)},f^{(m)},\hbar,l,u)$ \eqref{eq:mainDMR} admits a unique solution $(Y^{(m)}_t,Z^{(m)},K^{(m)},R^{(m)})\in S_G^\infty(0,T)\x\Hs_G^T(\R^d)\x\As_G^T\x C_B^T$. By Theorem 5.3 in \cite{QGBSDEbt} and Remark \ref{remark:QGZ}, the Q$G$-BSDE:
   \begin{equation} \label{eq:QGBSDEm}
   	\Yb^{(m)}_t=\xi^{(m)}+\int_t^Tf(s,Y^{(m)}_s,\Zb_s^{(m)})ds-\int_t^T\Zb^{(m)}_sdB_s-(\Kb^{(m)}_T-\Kb^{(m)}_t),~t\in[0,T],
   \end{equation}
 admits a unique solution $(\Yb^{(m)}_t,\Zb^{(m)}_t,\Kb^{(m)}_t)\in S_G^\infty(0,T)\x\Hs_G^T(\R^d)\x\As_G^T$ satisfying $(Y^{(m)},R^{(m)})=\widehat{\B\Sb\P}_l^u(\Yb^{(m)},\hbar)$ by Theorem \ref{theorem:BSPE&U}.
Meanwhile, following estimate \eqref{eq:fbound}, we have
\begin{equation}\label{eq:fmbound}
	\abv{f^{(m)}(t,\om,y,z)}\leq\ha_t(\om)+\frac{\gamma}{2}+\lm\abv{y}+\frac{3\gamma}{2}\abv{z}^2,~\forall m\in\N^*.
\end{equation}

\noindent{\bf Step 1. Uniform estimate.}
Set
\begin{align*}
	&\Dt:=\exp\crl{3p\gamma\sigmat^2\brak{\int_0^T\ha_sds+\frac{\gamma T}2+\lm T\sup_{s\in[0,T]}\max\edg{\abv{L_s(0)},\abv{U_s(0)}}}},\\
	&\Dt_0:=\exp\crl{3p\gamma\sigmat^2\sup_{s\in[0,T]}\max\edg{\abv{L_s(0)},\abv{U_s(0)}}}.
\end{align*}
Applying Lemma \ref{lemma:QGBSDEaprY} to Q$G$BSDE \eqref{eq:QGBSDEm} (taking $\hb_t:=\ha_t+\frac\gamma2+\lm\vert Y^{(n)}_t\vert$ and $\kappa:=3\gamma$) and then Proposition \ref{proposition:BSPSEapr}, it holds that for $t\in[0,T]$ and $p\geq1$,
\begin{align*}
	&\quad\Ebh\edg{\exp\crl{3p\gamma\sigmat^2\sup_{s\in[t,T]}\abv{\Yb^{(m)}_s}}}\\
	&\leq \Ebh\edg{\sup_{v\in[t,T]}\Ebh_v\edg{\exp\crl{3p\gamma\sigmat^2\brak{\abv{\xi}+\int_0^T\ha_sds+\frac{\gamma T}{2}+\lm (T-t)\sup_{s\in[t,T]}\abv{Y^{(m)}_s}}}}}\\
	&\leq \Ebh\Bigg[\sup_{v\in[t,T]}\Ebh_v\bigg[\exp\bigg\{3p\gamma\sigmat^2\bigg[\abv{\xi}+\int_0^T\ha_sds+\frac{\gamma T}{2}\\
	&~~~~~~~~~~~~~~~~~~~~~~~~~~~~~~~~~~~~~~~~~~+\lm (T-t)\bigg(\sup_{s\in[t,T]}\abv{\Yb^{(m)}_s}+\Cs\sup_{s\in[t,T]}\Ebh\edg{\abv{\Yb^{(m)}_s}}+C_0(t,T)\bigg)\bigg]\bigg\}\bigg]\Bigg]\\
    &\leq \Ebh\edg{\exp\crl{3p\gamma\sigmat^2\lm(T-t)\Cs\sup_{s\in[t,T]}\abv{\Yb^{(m)}_s}}}\brak{\Ebh\edg{\sup_{v\in[t,T]}\Ebh_v\edg{\exp\crl{6p\gamma\sigmat^2\brak{\lm (T-t)\sup_{s\in[t,T]}\abv{\Yb^{(m)}_s}}}}}}^\half\\
	&\quad\x\brak{\Ebh\edg{\sup_{s\in[t,T]}\Ebh_s\edg{\brak{\exp\crl{3p\gamma\sigmat^2\abv{\xi}}\Dt}^2}}}^\half\\
	&\leq \Ah_G\brak{\Ebh\edg{\brak{\exp\crl{3p\gamma\sigmat^2\abv{\xi}}\Dt}^4}}^\half\Ebh\edg{\exp\crl{6(2+\Cs)p\gamma\sigmat^2\lm (T-t)\sup_{s\in[t,T]}\abv{\Yb^{(m)}_s}}}.\notag
\end{align*}
 We select $h_3\in(0,T]$ such that $2(2+\Cs)\lm h_3<1$ and $T=N_h^3 h_3$, for some $N_h^3\in\N^*$ to get
\begin{align*}
	&\Ebh\edg{\exp\crl{3p\gamma\sigmat^2\sup_{s\in[T-h_3,T]}\abv{\Yb^{(m)}_s}}}\\
	&~~~~~~~~~~~\leq\brak{\Ah_G}^{\frac1{1-2(2+\Cs)\lm h_3}}\brak{\Ebh\edg{\brak{\exp\crl{24p\gamma\sigmat^2\abv{\xi}}}}}^{\frac1{4(1-2(2+\Cs)\lm h_3)}}\brak{\Ebh\edg{\Dt^8}}^{\frac1{4(1-2(2+\Cs)\lm h_3)}}.
\end{align*}
Then we get as $\Cs>1$,
\begin{align*}
	&\quad\Ebh\edg{\exp\crl{3p\gamma\sigmat^2\sup_{s\in[T-h_3,T]}\abv{Y^{(m)}_s}}}\leq \crl{\Ebh\edg{\exp\crl{3p\gamma\sigmat^2\sup_{s\in[T-h_3,T]}\abv{\Yb^{(m)}_s}}}}^{1+\Cs}\Dt_0\\
	&\leq(\Ah_G)^{\frac{1+\Cs}{1-2(2+\Cs)\lm h_3}}\brak{\Ebh\edg{\exp\crl{24p\gamma\sigmat^2\abv{\xi}}}}^{\frac{1+\Cs}{4(1-2(2+\Cs)\lm h_3)}}\brak{\Ebh\edg{\Dt^8}}^{\frac{1+\Cs}{4(1-2(2+\Cs)\lm h_3)}} \Dt_0:=\Dc(\xi,p).
\end{align*}
Next, iterating backwardly over the time intervals $[(k-1)h_3,kh_3]$ with terminal conditions $Y^{(m)}_{kh_3},~k=N_h^3-1,\cdots,1$, we have,
\begin{align*}
	&\quad\Ebh\edg{\exp\crl{3p\gamma\sigmat^2\sup_{s\in[(k-1)h_3,kh_3]}\abv{Y^{(m)}_s}}}\\
	&\leq(\Ah_G)^{\frac{1+\Cs}{1-2(2+\Cs)\lm h_3}}\brak{\Ebh\edg{\exp\crl{24p\gamma\sigmat^2\abv{Y^{(m)}_{kh_3}}}}}^{\frac{1+\Cs}{4(1-2(2+\Cs)\lm h_3)}}\brak{\Ebh\edg{\Dt^8}}^{\frac{1+\Cs}{4(1-2(2+\Cs)\lm h_3)}} \Dt_0\notag\\
	&\leq(\Ah_G)^{\frac{1+\Cs}{1-2(2+\Cs)\lm h_3}}\brak{\Dc(Y^{(n)}_{(k+1)h_3},8p)}^{\frac{1+\Cs}{4(1-2(2+\Cs)\lm h_3)}}\brak{\Ebh\edg{\Dt^8}}^{\frac{1+\Cs}{4(1-2(2+\Cs)\lm h_3)}} \Dt_0=:\Dc(Y^{(n)}_{kh_3},p),\notag
\end{align*}
whose boundedness is owing to Assumption \ref{assumption:HQUbdd}-($ii$). Then we obtain the uniform estimates:
\begin{align}
	&\quad\sup_{m\in\N^*}\Ebh\edg{\exp\crl{3p\gamma\sigmat^2\sup_{s\in[0,T]}\abv{Y^{(m)}_s}}}\leq\sup_{m\in\N^*}\Ebh\edg{\sup_{1\leq k\leq N_h^3}\exp\crl{3p\gamma\sigmat^2\sup_{s\in[(k-1)h_3,kh_3]}\abv{Y^{(m)}_s}}}\notag\\
	&\leq\sup_{m\in\N^*}\sum_{j=1}^{N_h^3}\Dc(Y^{(m)}_{kh_3},p)=:\Dch(\xi,p).\label{eq:unifesty}
\end{align}

\noindent{\bf Step 2. Convergence on a small time interval.}
Fix $m,q\geq1$ and $\theta\in(0,1)$, we define for $\phi=Y,~\Yb,~R^T-R,~\xi$,
\begin{align*}
	\dt\phi^{(m,q)}:=\frac{\phi^{(m+q)}-\theta\phi^{(m)}}{1-\theta},~\dtt\phi^{(m,q)}:=\frac{\phi^{(m)}-\theta\phi^{(m+q)}}{1-\theta},~\mbox{and}~~\dta\phi^{(m,q)}:=\abv{\dt\phi^{(m,q)}}+\abv{\dtt\phi^{(m,q)}}.
\end{align*}
Then the triple of processes $(\dt\Yb^{(m,q)},\dt\Zb^{(m,q)},\frac1{1-\theta}\Kb^{(m+q)})$ satisfies the following $G$-BSDE:
\begin{align*}
	\dt\Yb^{(m,q)}_t&=\dt\xi^{(m,q)}+\frac\theta{1-\theta}\brak{\Kb^{(m)}_T-\Kb^{(m)}_t}+\int_t^T\edg{\dt f^{(m,q)}\brak{s,\dt Y^{(m,q)}_s,\dt\Zb_s^{(m,q)}}+\dt \fb^{(m,q)}_s}ds\\
	&\quad-\int_t^T\dt\Zb^{(m,q)}_sdB_s-\frac1{1-\theta}\brak{\Kb^{(m+q)}_T-\Kb^{(m+q)}_t},
\end{align*}
where 
\begin{align*}
	&\dt\xi^{(m,q)}=\xi^{(m)}+\frac1{1-\theta}\brak{\xi^{(m+q)}-\xi^{(m)}}\leq\abv{\xi}+\frac1{1-\theta}\brak{\abv{\xi}-m}^+,\\ &\fb_t^{(m,q)}:=\frac{\fb^{(m+q)}_t-\theta\fb^{(m)}_t}{1-\theta}-\fb_t\leq\fb_t^{(m)}-\fb_t+\frac1{1-\theta}\brak{\abv{\fb_t}-m}^+\leq\frac2{1-\theta}\brak{\abv{\fb_t}-m}^+,
\end{align*}
and by Assumption \ref{assumption:HQUbdd}-($iii$)($v$) with estimate \eqref{eq:fbound},
\begin{align*}
	&\quad\dt f^{(m,q)}(t,\dt Y^{(m,q)}_t,z)\\
	&\leq\lm\abv{\dt Y^{(m,q)}_t}+\lm\abv{Y^{(m)}_t}+\frac1{1-\theta}\edg{f(t,Y^{(m)}_t,(1-\theta)z+\theta\Zb^{(m)}_t)-\theta f(t,Y^{(m)}_t,\Zb^{(m)}_t)}\\
	&\leq \lm\abv{\dt Y^{(m,q)}_t}+\lm\abv{Y^{(m)}_t}+f(t,Y^{(m)}_t,z)\leq \a_t+\frac\gamma2+2\lm\abv{Y^{(m)}_t}+\lm\abv{\dt Y^{(m,q)}_t}+\frac{3\gamma}{2}\abv{z}^2.
\end{align*}
Set
\begin{align*}
	\Phi^{(m,q)}&:=\int_0^T\a_sds+\frac{\gamma T}{2}+2\lm T\brak{\sup_{s\in[0,T]}\abv{Y^{(m)}_s}+\sup_{s\in[0,T]}\abv{Y^{(m+q)}_s}},\\
	\Gamma^{(m,q)}&:=\abv{\xi}+\int_0^T\a_sds+\frac{\gamma T}{2}+\lm T\brak{\sup_{s\in[0,T]}\abv{Y^{(m)}_s}+\sup_{s\in[0,T]}\abv{Y^{(m+q)}_s}},\\
	\mbox{and}~~\rho(\theta,m)&:=\frac1{1-\theta}\brak{\abv{\xi}-m}^++\frac2{1-\theta}\int_0^T\brak{\abv{\fb_s}-m}^+ds.
\end{align*}
Note that in the spirit of Assumption \ref{assumption:HQUbdd}-($ii$) and estimate \eqref{eq:unifesty}, $\Phi^{(m,q)}$ and $\Gamma^{(m,q)}$ admit finite exponential moments of arbitrary order, $\forall m,q\in\N^*$.
Applying Lemma \ref{lemma:QGBSDEaprY}-($ii$) (taking $\hb_t:=\ha_t+\frac\gamma2+2\lm\abv{Y^{(m)}_t}+\lm\abv{\dt Y^{(m,q)}_t}$ and $\kappa:=3\gamma$), we get by symmetry,
\begin{align*}
	&\max\brak{\exp\crl{3p\gamma\sigmat^2\brak{\dt\Yb^{(m,q)}_t}^+},\exp\crl{3p\gamma\sigmat^2\brak{\dtt\Yb^{(m,q)}_t}^+}}\\
	&~~~~~~~~~~~~~~~~~~~~~~~~~~\leq\Ebh_t\edg{\exp\crl{3p\gamma\sigmat^2\edg{\abv{\xi}+\rho(\theta,m)+\Phi^{(m,q)}+\lm(T-t)\sup_{s\in[t,T]}\abv{\dta Y^{(m,q)}_s}}}}.
\end{align*}
Considering the fact that
\begin{align*}
	\brak{\dt \Yb^{(m,q)}}^-\leq \brak{\dtt\Yb^{(m,q)}}^++\abv{\Yb^{(m)}}+\abv{\Yb^{(m+q)}}~\mbox{and}~~\brak{\dtt\Yb^{(m,q)}}^-\leq\brak{\dt\Yb^{(m,q)}}^++\abv{\Yb^{(m)}}+\abv{\Yb^{(m+q)}},
\end{align*}
we apply estimate \eqref{eq:fmbound}, Lemma \ref{lemma:QGBSDEaprY} and Corollary \ref{corollary:Doob} to deduce that
\begin{align}
	&\quad\Ebh\edg{\exp\crl{3p\gamma\sigmat^2\sup_{s\in[t,T]}\dta\Yb^{(m,q)}_s}}\notag\\
	&\leq\Ebh\edg{\brak{\exp\crl{3p\gamma\sigmat^2\sup_{s\in[t,T]}\max\brak{\abv{\dt\Yb^{(m,q)}_s},\abv{\dtt\Yb^{(m,q)}_s}}}}^2}\notag\\
	&\leq \Ebh\edg{\brak{\exp\crl{3p\gamma\sigmat^2\edg{\brak{\dt \Yb^{(m,q)}_t}^++\brak{\dtt \Yb^{(m,q)}_t}^++\abv{\Yb^{(m)}_t}+\abv{\Yb^{(m+q)}_t}}}}^2}\notag\\
	&\leq\Bigg\{\Ebh\Bigg[\sup_{v\in[t,T]}\brak{\Ebh_v\edg{\exp\crl{3p\gamma\sigmat^2\edg{\abv{\xi}+\rho(\theta,m)+\Phi^{(m,q)}+\lm(T-t)\sup_{s\in[t,T]}\dta Y^{(m,q)}_s}}}}^8\Bigg]\notag\\
	&\quad\quad\x\Ebh\edg{\sup_{v\in[t,T]}\brak{\Ebh_v\edg{\exp\crl{3p\gamma\sigmat^2\Gamma^{(m,q)}}}}^8}\Bigg\}^\half\label{eq:Ynkexp-0}
	\\
	&\leq\Ah_G\crl{\Ebh\edg{\exp\crl{48p\gamma\sigmat^2\edg{\abv{\xi}+\rho(\theta,m)+\Phi^{(m,q)}+\lm(T-t)\sup_{s\in[t,T]}\dta Y^{(m,q)}_s}}}\Ebh\edg{\exp\crl{48p\gamma\sigmat^2\Gamma^{(m,q)}}}}^\half.\notag
\end{align}
From Lemma \ref{lemma:QGBSDEaprY}, we have for each $m\in\N^*$,
\begin{align}
	&\quad\Ebh\edg{\exp\crl{3p\gamma\sigmat^2\sup_{s\in[t,T]}\abv{\Yb_s^{(m)}}}}\notag\\
	&\leq \Ebh\edg{\sup_{s\in[t,T]}\Ebh_s\edg{\exp\crl{3p\gamma\sigmat^2\brak{\abv{\xi^{(m)}}+\int_0^T\ha_vdv+\frac{\gamma T}{2}+\lm T\sup_{v\in[t,T]}\abv{Y^{(m)}_v}}}}}\notag\\
	&\leq \Ah_G\Ebh\edg{\exp\crl{6p\gamma\sigmat^2\brak{\int_0^T\ha_sds+\frac{\gamma T}{2}+(1+\lm T)\sup_{s\in[t,T]}\abv{Y^{(m)}_s}}}}\notag\\
	&\leq \Ah_G\brak{\Ebh\edg{\exp\crl{12p\gamma\sigmat^2\brak{\int_0^T\ha_sds+\frac{\gamma T}{2}}}}}^\half\brak{\Dch(\xi,\lceil4(1+\lm T)\rceil p)}^\half.\label{eq:Ybexpn}
\end{align}
Then, following Lemma \ref{lemma:Rtheta}, we could set
\begin{align}
	&\quad C_0'(t,T,\Yb^{(m)},\Yb^{(m+k)})\notag\\
	&\leq\frac{2}{\cu}\sup_{s\in[0,T]}\max\brak{\abv{l_s},\abv{u_s}}+\frac{\cb}{\cu}\sup_{s\in[0,T]}\max\brak{\abv{L_s(0)},\abv{U_s(0)}}\label{eq:Cs0star}\\
	&\quad+\frac{4\cb(\cu+\cb)\Ah_G}{3\cu^2\gamma\sigmat^2}\brak{\Ebh\edg{\exp\crl{12\gamma\sigmat^2\brak{\int_0^T\ha_sds+\frac{\gamma T}{2}}}}}^\half\brak{\Dch(\xi,\lceil4(1+\lm T)\rceil )}^\half=:\Cs_0^*,\notag
\end{align}
of which the boundedness is guaranteed by Proposition \ref{proposition:ULcon}-($i$) and estimate \eqref{eq:unifesty} for any $(m,q)\in\N^*\x\N^*$. 
Applying Lemma \ref{lemma:Rtheta} and estimate \eqref{eq:Cs0star} in the first inequality, we have 
\begin{align*}
	&\quad\Ebh\edg{\exp\crl{3p\gamma\sigmat^2\sup_{s\in[t,T]}\dta Y^{(m,q)}_s}}\\
	&\leq\Ebh\edg{\exp\crl{3p\gamma\sigmat^2\brak{2\Cs_0^*+\sup_{s\in[t,T]}\dta\Yb^{(m,q)}_s+2\Cs\sup_{s\in[t,T]}\Ebh\edg{\dta\Yb^{(m,q)}_s}}}}\\
	&\leq \Ebh\edg{\exp\crl{6(1+2\Cs)p\gamma\sigmat^2\brak{\sup_{s\in[t,T]}\dta \Yb^{(m,q)}_s+\Cs_0^*}}}\\
	&\leq \Ah_G\brak{\Ebh\edg{\exp\crl{96(1+2\Cs)p\gamma\sigmat^2\edg{\abv{\xi}+\rho(\theta,m)+\Cs_0^*+\Phi^{(m,q)}+\lm(T-t)\sup_{s\in[t,T]}\dta Y^{(m,q)}_s}}}}^\half\\
	&\quad\x\brak{\Ebh\edg{\exp\crl{96(1+2\Cs)p\gamma\sigmat^2\Gamma^{(m,q)}}}}^\half\\
	&\leq A_G^*(p)\brak{\Ebh\edg{\exp\crl{384(1+2\Cs)p\gamma\sigmat^2\rho(\theta,m)}}}^{\frac18}\\
	&\quad\x\brak{\Ebh\edg{\exp\crl{192(1+2\Cs)p\gamma\sigmat^2\lm(T-t)}\sup_{s\in[t,T]}\dta Y^{(m,q)}_s}}^\frac14,
\end{align*}
where the second inequality is deduced following the procedure of estimate \eqref{eq:Ynkexp-0} and
\begin{align*}
	&\quad A^*_G(p)\\
	&:= \Ah_G\brak{\Ebh\edg{\exp\crl{384(1+2\Cs)p\gamma\sigmat^2\edg{\abv{\xi}+\Cs_0^*+\Phi^{(m,q)}}}}}^{\frac18}\brak{\Ebh\edg{\exp\crl{96(1+2\Cs)p\gamma\sigmat^2\Gamma^{(m,q)}}}}^\half\\
	&<+\infty.
\end{align*}
Choosing $h_4>0$ such that $64(1+2\Cs)\lm h_4<1$ and $T=N_h^4h_4$ for some $N_h^4\in\N^*$, we have
\begin{align}
	&\quad\Ebh\edg{\exp\crl{3p\gamma\sigmat^2\sup_{s\in[T-h_4,T]}\dta Y^{(m,q)}_s}}\label{eq:Ynkexp-1}\\
	&\leq \brak{A^*_G(p)}^{\frac1{1-16(1+2\Cs)\lm h_4}}\notag\\
	&\quad\x\brak{\Ebh\edg{\exp\crl{\frac{384(1+2\Cs)p\gamma\sigmat^2}{1-\theta}\brak{\brak{\abv{\xi}-m}^++2\int_0^T\brak{\abv{\fb_s}-m}^+ds}}}}^{\frac1{8(1-16(1+2\Cs)\lm h_4)}}.\notag
\end{align}
Since by Assumption \ref{assumption:HQUbdd}-($ii$)($iii$),
\begin{align*}
	L_G^1(\Om_T)\ni\exp\crl{\frac{384(1+2\Cs)p\gamma\sigmat^2}{1-\theta}\brak{\brak{\abv{\xi}-m}^++2\int_0^T\brak{\abv{\fb_s}-m}^+ds}}\downarrow1,~\mbox{as}~n\to\infty,
\end{align*}
following monotone convergence theorem, for each $p\geq1$ and any $\theta\in(0,1)$, we have
\begin{equation}\label{eq:Ynkexp-3}
	\limsup_{m\to\infty}\sup_{q\in\N^*}\Ebh\edg{\exp\crl{3p\gamma\sigmat^2\sup_{s\in[T-h_4,T]}\dta Y^{(m,q)}_s}}\leq\brak{A^*_G(p)}^{\frac1{1-16(1+2\Cs)\lm h_4}}<+\infty.
\end{equation}
Because
$
	Y^{(m+q)}-Y^{(m)}=(1-\theta)\brak{\dt Y^{(m,q)}-Y^{(m)}}$,
it holds that for any $k\in\N^*$,
\begin{align*}
	&\quad\limsup_{m\to\infty}\sup_{q\in\N^*}\Ebh\edg{\sup_{s\in[T-h_4,T]}\abv{ Y^{(m+q)}_s-Y^{(m)}_s}^k}\\
	&\leq 2^{k-1}(1-\theta)^k\brak{\limsup_{m\to\infty}\sup_{q\in\N^*}\Ebh\edg{\sup_{s\in[T-h_4,T]}\abv{\dt Y^{(m,q)}_s}^k}+\sup_{m\in\N^*}\Ebh\edg{\sup_{s\in[T-h_4,T]}\abv{Y^{(m)}_s}^k}}.
\end{align*}
Letting $\theta\to1$, we can assert by the bounds \eqref{eq:unifesty}-\eqref{eq:Ynkexp-3} that $\crl{Y^{(m)}}_{m\in\N}$ is a Cauchy sequence in $S_G^q(T-h_4,T),~\forall q\geq1$, which indicates that there exists a continuous process $\Yc^{N_h^4}$ such that
\begin{equation*}
	\lim_{m\to\infty}\Ebh\edg{\sup_{s\in[T-h_4,T]}\abv{Y^{(m)}_s-\Yc_s^{N_h^4}}^n}=0,~\forall n\geq1.
\end{equation*}

\noindent{\bf Step 3. Convergence on the whole time horizon.}

Following Step 2, we can obtain a sequence $\Yc^{k},~1\leq k\leq N_h^4$ such that
\begin{equation}\label{eq:Ycjcvg}
	\lim_{m\to\infty}\Ebh\edg{\sup_{s\in[(k-1)h_4,kh_4]}\abv{Y^{(m)}_s-\Yc^{j}_s}^n}=0,~\forall n\geq1.
\end{equation}
On each time intervals $[(k-1)h_4,kh_4],~k=1,\cdots,N^4_h-1$, 
we replace $\Gamma^{(m,q)}$ and $\rho(\theta,m)$ by
\begin{align*}
	&\Gamma^{(m,q)}_k:=\max\brak{\abv{Y^{(m)}_{kh_4}},\abv{Y^{(m+q)}_{kh_4}}}+\int_0^T\a_sds+\frac{\gamma T}{2}+\lm T\brak{\sup_{s\in[0,T]}\abv{Y^{(m)}_s}+\sup_{s\in[0,T]}\abv{Y^{(m+q)}_s}},\\
	&\rho(\theta,m,q,k):=\frac1{1-\theta}\brak{\abv{Y^{(m)}_{kh_4}-Y^{(m+q)}_{kh_4}}+2\int_0^T\brak{\abv{\fb_s}-m}^+ds}.
\end{align*}
Hereafter, $\dt Y^{(m,q)}_{kh_4}$ will take the role of $\dt \xi^{(m,q)}$ on $[T-h_4,T]$. Set
\begin{align}
	&\quad\Ab_G^{*}(p)\label{eq:Ynkexp-7}\\
	&:=\Ah_G\sup_{m,q\in\N^*,1\leq k\leq N_h^4}\brak{\Ebh\edg{\exp\crl{384(1+2\Cs)p\gamma\sigmat^2\edg{\max\brak{\abv{Y^{(m)}_{kh_4}},\abv{Y^{(m+q)}_{kh_4}}}+\Cs_0^*+\Phi^{(m,q)}}}}}^{\frac18}\notag\\
	&~\quad\x\sup_{m,q\in\N^*,1\leq k\leq N_h^4}\brak{\Ebh\edg{\exp\crl{96(1+2\Cs)p\gamma\sigmat^2\Gamma^{(m,q)}_k}}}^\half<\infty,\notag
\end{align}
which is uniformly bounded using estimate \eqref{eq:unifesty}. Thereby, similarly to \eqref{eq:Ynkexp-1}, we have
\begin{align}\label{eq:Ynkexp-9}
	&\quad\lim_{m\to\infty}\sup_{q\in\N^*}\Ebh\edg{\exp\crl{3p\gamma\sigmat^2\sup_{s\in[(k-1)h_4,kh_4]}\dta Y^{(m,q)}_s}}\\
	&\leq \brak{\Ab_G^*(p)}^{\frac1{1-16(1+2\Cs)\lm h_4}}\lim_{m\to\infty}\sup_{q\in\N^*}\brak{\Ebh\edg{\exp\crl{384(1+2\Cs)p\gamma\sigmat^2\rho(\theta,m,q,k)}}}^{\frac1{8(1-16(1+2\Cs)\lm h_4)}}\notag\\
	&\leq \brak{\Ab_G^*(p)}^{\frac1{1-16(1+2\Cs)\lm h_4}}\lim_{m\to\infty}\sup_{q\in\N^*}\brak{\Ebh\edg{\exp\crl{768(1+2\Cs)p\gamma\sigmat^2\brak{\abv{Y^{(m)}_{kh_4}}+\dta Y^{(m,q)}_{kh_4}}}}}^{\frac1{16(1-16(1+2\Cs)\lm h_4)}},\notag
\end{align}
which is uniformly bounded by estimates for \eqref{eq:unifesty} and
$$\lim_{m\to\infty}\sup_{q\in\N^*}\Ebh\edg{\exp\crl{3p\gamma\sigmat^2\sup_{s\in[(k'-1)h_4,~k'h_4]}\dta Y^{(m,q)}_s}},~k+1\leq k'\leq N_h^4.$$
Then assertion \eqref{eq:Ycjcvg} follows, and we put
$$
	Y_t:=\sum_{k=1}^{N_h^4}\Yc^{k}_t\1_{[(k-1)h_4,kh_4)}(t)+\xi\1_{\crl{T}}(t),~t\in[0,T].$$
By the convergence on each small intervals, it holds that $Y$ has a quasi-continuous version, which is still denoted by $Y$. Then it holds that
\begin{align}
	\lim_{m\to\infty}\Ebh\edg{\sup_{t\in[0,T]}\abv{Y^{(m)}_t-Y_t}^n}\leq\sum_{k=1}^{N_h^4}\lim_{m\to\infty}\Ebh\edg{\sup_{t\in[(k-1)h_4,kh_4]}\abv{Y^{(m)}_t-\Yc_t^{k}}^n}=0.\label{eq:Ytotalcvg}
\end{align} 
Moreover, with the aid of estimate \eqref{eq:unifesty} and Fatou's Lemma (Theorem \ref{theorem:mcvg}), $Y\in\Es_G^T$ and 
\begin{equation}\label{eq:Yexpest}
	\Ebh\edg{\exp\crl{3p\gamma\sigmat^2\sup_{t\in[0,T]}\abv{Y_t}}}\leq \Dch(\xi,p),~\forall p\geq1.
\end{equation} 

\noindent{\bf Step 4. Verification of the solution.}

By Theorem 3.9 in \cite{QGBSDEubt}, the following Q$G$-BSDE:
\begin{equation*}
	\Yb_t=\xi+\int_t^Tf(s,Y_s,\Zb_s)ds-\int_t^T\Zb_sdB_s-(\Kb_T-\Kb_t),~t\in[0,T],
\end{equation*} 
admits a unique solution $(\Yb,\Zb,\Kb)\in\Es_G^T\x\Hs_G^T(\R^d)\x\As_G^T$, where the bounds for $\\ \Ebh[\exp\{3p\gamma\sigmat^2\sup_{t\in[0,T]}\abv{\Yb_t}\}],~p\geq1$ can be obtained via Lemma \ref{lemma:QGBSDEaprY}.
From Theorem \ref{theorem:BSPE&U}, we can acquire uniquely $C_G^1(0,T)\x C_B^T\ni(\Yt,\Rt)=\widehat{\B\Sb\P}_l^u(\Yb,\hbar)$. Therefore the quadruple $(\Yt,\Zb,\Kb,\Rt)$ satisfies the DMR$G$-BSDE:
\begin{equation*}
	\left\{
	\begin{lgathered}
		\Yt_t=\xi+\int_t^Tf(s,Y_s,\Zb_s)ds-\int_t^T\Zb_sdB_s-(\Kb_T-\Kb_t)+(\Rt_T-\Rt_t),~t\in[0,T];\\
		l_t\leq\Ebh\edg{\hbar(t,\Yt_t)}\leq u_t,~\forall t\in[0,T];\\
		\int_{s}^{t}\brak{\Ebh\edg{\hbar(v,\Yt_v)}-l_v}dR_v\leq0,~\mbox{and}~~\int_{s}^{t}\brak{\Ebh\edg{\hbar(v,\Yt_v)}-u_v}dR_v\leq0,~ \forall0\leq s\leq t\leq T.
	\end{lgathered}
	\right.
\end{equation*}
Following Proposition \ref{proposition:BSPSEapr}, $\Yt\in\Es_G^T$.
It remains to check that $Y=\Yt$.
 For each $m\in\N^*$, 
\begin{align}
	\Ebh\edg{\sup_{t\in[0,T]}\abv{\Yt_t-Y_t}}&\leq\Ebh\edg{\sup_{t\in[0,T]}\abv{\Yt_t-Y^{(m)}_t}}+\Ebh\edg{\sup_{t\in[0,T]}\abv{Y^{(m)}_t-Y_t}}\label{eq:YYnsup}\\
	&\leq (1+\Cs)\Ebh\edg{\sup_{t\in[0,T]}\abv{\Yb_t-\Yb^{(m)}_t}}+\Ebh\edg{\sup_{t\in[0,T]}\abv{Y^{(m)}_t-Y_t}},\notag
\end{align} 
where the last inequality is attributed to Proposition \ref{proposition:BSPSEaprDf}.

Since by estimate \eqref{eq:Ybexpn} we can get a uniform bound for $\sup_{m\in\N^*}\Ebh[\exp\{3p\gamma\sigmat^2\sup_{t\in[0,T]}\vert\Yb^{(m)}\vert\}]$. Following a $\theta$-method for ($Y^{(m)},Y,\Yb^{(m)},\Yb$) on $[0,T]$ similar to that for ($Y^{(m)},Y^{(m+q)},\Yb^{(m)},\Yb^{(m+q)}$) in Step 2, we get the convergence result to zero of the first term in estimate \eqref{eq:YYnsup}. Then applying the convergence result of $\brak{Y^{(m)}}_{m\in\N^*}$ to $Y$, we can complete the proof for existence. Moreover, the uniqueness issue can be obtained analogously to Lemma 4.4 in \cite{GLX}, so we omit it.
\qed 
\section{Multi-dimensional doubly mean-reflected $G$-BSDEs}
In this section, we consider a group of multi-dimensional DMR$G$-BSDEs: for each $1\leq j\leq n$,
\begin{equation}\label{eq:mdDMRGBSDE}
	\left\{
	\begin{lgathered}
		Y^j_t=\xi^j+\int_t^Tf^j(s,Y_s,Z^j_s)ds-\int_t^TZ^j_sdB_s-(K^j_T-K^j_t)+(R^j_T-R^j_t),~t\in[0,T];\\
		l^j_t\leq\Ebh\edg{\hbar^j(t,Y^j_t)}\leq u^j_t,~\forall t\in[0,T];\\
		\int_{s}^{t}\brak{\Ebh\edg{\hbar^j(v,Y^j_v)}-l^j_v}dR^j_v\leq0,~\mbox{and}~~\int_{s}^{t}\brak{\Ebh\edg{\hbar^j(v,Y^j_v)}-u^j_v}dR^j_v\leq0,~\forall0\leq s\leq t\leq T;
	\end{lgathered}
	\right.
\end{equation}
where the generators 
$$\fv(t,\om,y,z):=(f^1(t,\om,y,z^1),\cdots,f^n(t,\om,y,z^n))^\top:[0,T]\x\Om_T\x\R^n\x\R^{n\x d}\to\R^n,$$
with $y^j$ and $z^j$ denoting the $j$-th component of $y$ and the $j$-th row of $z$ for any $(y,z)\in\R^n\x\R^{n\x d}$, respectively. 
For the ease of symbolization, we denote 
$$\xi:=(\xi^1,\cdots,\xi^n)^\top,~\hv(t,y):=(\hbar^1(t,y^1),\cdots,\hbar^n(t,y^n))^\top, ~\lv:=(l^1,\cdots,l^n)^\top,~\mbox{and}~~\uv:=(u^1,\cdots,u^n)^\top.$$
Suppose the following assumptions are imposed on the data $(\xi,\fv,\hv,\lv,\uv)$. Firstly, a universal hypothesis on reflections and terminal conditions is stated
\begin{Assumption}\label{assumption:mdhrxi}
	\begin{enumerate}[($i$)]
		\item For each $1\leq j\leq n$, $(\hbar^j,l^j,u^j)$ satisfies Assumption \ref{assumption:HR} with the same constants $0<\cu\leq \cb$.
		\item $\xi\in L_G^1(\Om_T;\R^n)$ and $l^j_T\leq \Ebh\edg{\hbar^j(T,\xi^j)}\leq u^j_T,~\forall1\leq j\leq n$.
	\end{enumerate}
\end{Assumption}
Then, the Lipschitz generators are assumed in the form:
\begin{Assumption}\label{assumption:mdHLip}
	Fix some constants $\b>1$ and $L>0$.
	\begin{enumerate}[($i$)]
		\item $\xi\in L_G^\b(\Om_T;\R^n)$.
		\item For each $1\leq j\leq n$ and any $(t,\om)\in[0,T]\x\Om_T$, $f^j(t,\om,0,\0)\in M_G^\b(0,T)$.
		\item For each $(t,\om)\in[0,T]\x\Om_T$ and any $(y,y',z,z')\in\R^{2n+2n\x d}$,
		\begin{equation*}
			\abv{\fv(t,\om,y,z)-\fv(t,\om,y',z')}\leq L\brak{\abv{y-y'}+\abv{z-z'}}.
		\end{equation*}
	\end{enumerate}
\end{Assumption}
Next, we consider the quadratic generators with bounded and unbounded terminal conditions respectively.
\begin{Assumption}\label{assumption:MDqf}
	Let $\lm,\gamma\in\R^+$.
	\begin{enumerate}[($i$)]
		\item For any $(t,\om)\in[0,T]\x\Om_T$ and any  $(y,y',z,z')\in\R^{2n+2n\x d}$, 
		$$\abv{\fv(t,\om,y,z)-\fv(t,\om,y',z')}\leq\lm\abv{y-y'}+\gamma(1+\abv{z}+\abv{z'})\abv{z-z'}.$$
		\item There is a modulus of continuity $w:[0,\infty)\to[0,\infty)$ such that
		$$\abv{\fv(t,\om,y,z)-\fv(t',\om',y,z)}\leq w(\abv{t-t'}+\norm{\om-\om'}),$$
		for each $(y,z)\in\R^{n+n\x d}$ and any  $(t,t',\om,\om')\in[0,T]\x[0,T]\x\Om_T\x\Om_T$.
	\end{enumerate}
\end{Assumption}
\begin{Assumption}\label{assumption:mdHQbdd}Let $M_0\in\R^+$.
	\begin{enumerate}[($i$)]
		\item $\xi\in L_G^\infty(\Om_T;\R^n)$.
		\item There exists some constant $M_0>0$ such that for each $(t,\om)\in[0,T]\x\Om_T$,
		$$\abv{\xi}+\int_0^T\abv{\fv(s,\om,\0,\0)}^2ds\leq M_0.$$
	\end{enumerate}
\end{Assumption}
\begin{Assumption}\label{assumption:mdHQUbdd}
	Let $\ha\in M_G^1(0,T)$ be nonnegative.
	\begin{enumerate}[($i$)]
		\item Both $\xi$ and $\int_0^T\ha_tdt$ have finite exponential moments of arbitrary order, i.e.,
		$$\Ebh\edg{\exp\crl{p\abv{\xi}+p\int_0^T\ha_tdt}}<\infty,~\forall p\geq1.$$
		\item For each $1\leq j\leq n$ and any $(t,\om,y)\in[0,T]\x\Om_T\x\R^n$, $f^j(t,\om,y,\cdot)$ is convex.
	\end{enumerate}
\end{Assumption}
\begin{Theorem}
	Let Assumption \ref{assumption:mdhrxi} hold true.
	\begin{enumerate}[($i$)]
		\item If Assumption \ref{assumption:mdHLip} holds, then the multi-dimensional DMR$G$-BSDE ($\xi,\fv,\hv,\lv,\uv$) \eqref{eq:mdDMRGBSDE} admits a unique solution $(Y,Z,K,R)\in\Sc_G^\a(0,T;\R^n)\x C_B^T(\R^n)$, for any $1\leq \a<\b$.
		\item If Assumptions \ref{assumption:MDqf}-\ref{assumption:mdHQbdd} holds, then the multi-dimensional DMR$G$-BSDE ($\xi,\fv,\hv,\lv,\uv$) \eqref{eq:mdDMRGBSDE} admits a unique solution $(Y,Z,K,R)\in \S_G^\infty(0,T;\R^n)\x\Hs_G^T(\R^{n\x d})\x\As_G^T(\R^n)\x C_B^T(\R^n)$.
	\end{enumerate}
\end{Theorem}
\noindent{\bf\underline{Proof.}}
($i$) Taking $V^i\in M_G^\b(0,T;\R^n),~i=1,2$,
we consider the following one-dimensional DMR$G$-BSDE with the corresponding reflection and Skorokhod condition given by \eqref{eq:mdDMRGBSDE}:
\begin{equation}\label{eq:mdtest}
	Y^{i,j}_t=\xi^j+\int_t^Tf^j(s,V^i_s,Z^{i,j}_s)ds-\int_t^TZ^{i,j}_sdB_s-(K^{i,j}_T-K^{i,j}_t)+(R^{i,j}_T-R^{i,j}_t),~t\in[0,T],
\end{equation}
which admits a unique solution $(Y^{L,i,j},Z^{L,i,j},K^{L,i,j},R^{L,i,j})\in \Sc_G^\a(0,T)\x C_B^T$, $1\leq \a<\b$ by Lemmas 3.2-3.3 in \cite{LiumdGBSDE} and Theorem \ref{theorem:LipDMR}. Similarly to the proof of Theorem \ref{theorem:LipDMR}, we have by Proposition \ref{prop:1daprDY},
\begin{align*}
	\norm{Y^{L,1,j}-Y^{L,2,j}}^{\b}_{M_G^\b(t,T)}\leq 2^{\b-1}(1+\Cs^\b)C^LL^\b(T-t)^\b\norm{U^1-U^2}_{M_G^{\b}(t,T;\R^n)}^\b,
\end{align*}
and then
\begin{align*}
	\norm{Y^{L,1}-Y^{L,2}}^{\b}_{M_G^\b(t,T;\R^n)}\leq 2^{\b-1}n^{(\bh-1)\vee0}(1+\Cs^\b)C^LL^\b(T-t)^\b\norm{U^1-U^2}_{M_G^{\b}(t,T;\R^n)}^\b.
\end{align*}
Therefore, the contraction mapping and iteration process can be realized similarly to Lemma 3.4 and Theorem 3.1 in \cite{LiumdGBSDE}.

\noindent($ii$) Taking $V^i\in S_G^\infty(0,T;\R^n), i=1,2$ then the one-dimensional DMR$G$-BSDE \eqref{eq:mdtest}
admit a unique solution $(Y^{Q,i,j},Z^{Q,i,j},K^{Q,i,j},R^{Q,i,j})\in S_G^\infty(0,T)\x\Hs_G^T(\R^d)\x\As_G^T\x C_B^T$ from Theorem \ref{theorem:btv-2}. Then we compute by a similar $G$-Girsanov transformation and a priori estimates of backward Skorokhod problem with sublinear expectation to get that
\begin{align*}
	\norm{Y^{Q,1}-Y^{Q,2}}_{S_G^\infty(t,T;\R^n)}\leq n(1+\Cs)\lm(T-t)\norm{U^1-U^2}_{S_G^\infty(t,T;\R^n)},
\end{align*} 
and complete this assertion analogously in the sense of a fixed-point argumentation.
\qed

\begin{Theorem}
	Let Assumptions \ref{assumption:mdhrxi}-\ref{assumption:MDqf}-\ref{assumption:mdHQUbdd} hold. Then the multi-dimensional DMR$G$-BSDE ($\xi,\fv,\hv,\lv,\uv$) \eqref{eq:mdDMRGBSDE} admits a unique solution $(Y,Z,K,R)\in \Es_G^T(\R^n)\x\Hs_G^T(\R^{n\x d})\x\As_G^T(\R^n)\x C_B^T(\R^n)$.
\end{Theorem}
\noindent{\bf\underline{Proof.}}
\noindent{\bf Iteration.} Taking $V\in \Es_G^T(\R^n)$, the following quadratic $G$-BSDE:
\begin{equation*}
	\Yb^j_t=\xi^j+\int_t^Tf^j(s,V_s,\Zb_s^j)ds-\int_t^T\Zb_s^jdB_s-(\Kb^j_T-\Kb^j_t),~t\in[0,T],
\end{equation*}
admits a unique solution $(\Yb^j,\Zb^j,\Kb^j)\in \Es_G^T\x\Hs_G^T(\R^d)\x\As_G^T$ by Lemma 4.1 in \cite{QGBSDEubt}. From Theorem \ref{theorem:BSPE&U}, there exists $C_G^1(0,T)\x C_B^T\ni(Y^j,R^j)=\widehat{\B\Sb\P}_{l^j}^{u^j}(\Yb^j,\hbar^j)$ and $Y^j\in \Es_G^T$ from Proposition \ref{proposition:BSPSEapr}. Combining with Theorem \ref{theorem:main},  the quadruple $(Y^j,\Zb^j,\Kb^j,R^j)$ is the unique solution to DMR$G$-BSDE:
\begin{equation*}
	Y^j_t=\xi^j+\int_t^Tf^j(s,V_s,Z^j_s)ds-\int_t^TZ_s^jdB_s-(K^j_T-K^j_t)+(R^j_T-R^j_t),~t\in[0,T].
\end{equation*} 
Then we can apply an iterative method by setting $Y^{[0]}\equiv\0$ and solving the following  DMR$G$-BSDEs on $[0,T]$:
\begin{equation*}
	Y^{[m],j}_t=\xi^j+\int_t^Tf^j(s,Y^{[m-1]}_s,Z^{[m],j}_s)ds-\int_tZ^{[m],j}_sdB_s-(K^{[m],j}_T-K^{[m],j}_t)+(R^{[m],j}_T-R^{[m],j}_t),
\end{equation*}
and denote by
\begin{equation*}
	\Yb^{[m],j}_t=\xi^j+\int_t^Tf^j(s,Y^{[m-1]}_s,Z^{[m],j}_s)ds-\int_tZ^{[m],j}_sdB_s-(K^{[m],j}_T-K^{[m],j}_t).
\end{equation*}
For each $m\in\N^*$, $(Y^{[m]},Z^{[m]},K^{[m]},R^{[m]})\in\Es_G^T(\R^n)\x\Hs_G^T(\R^{n\x d})\x\As_G^T(\R^n)\x C_B^T(\R^n)$.
Due to Assumption \ref{assumption:mdHQUbdd}, we can get similarly to \eqref{eq:fbound} that
\begin{equation}
	\abv{f^j(t,V_t,z^j)}\leq \ha_t+\frac\gamma2+\lm\abv{V_t}+\frac{3\gamma}{2}\abv{z^j}^2.
\end{equation}
\noindent{\bf Uniform estimate.}
From assertion ($i$) of Lemma \ref{lemma:QGBSDEaprY} (taking $\hb_t=\ha_t+\frac\gamma2+\lm\abv{Y_t^{[m-1]}}$ and $\kappa=3\gamma$), we have that for $p\geq1$ and $1\leq j\leq n$,
\begin{align*}
	\exp\crl{3p\gamma\sigmat^2\abv{\Yb^{[m],j}_t}}\leq \Ebh_t\edg{\exp\crl{3p\gamma\sigmat^2\brak{\abv{\xi}+\int_t^T\brak{\ha_s+\frac\gamma2+\lm\abv{Y^{[m-1]}_s}}ds}}},~m\geq1.
\end{align*}
Recalling Proposition \ref{proposition:BSPSEapr}, we set $$\Dt_1(p):=\exp\crl{3p\gamma\sigmat^2\sum_{j=1}^n\sup_{s\in[0,T]}\max\edg{\abv{L^j_s(0)},\abv{U^j_s(0)}}},$$
and apply Jensen's inequality to get
\begin{align*}
	\exp\crl{3p\gamma\sigmat^2\abv{Y^{[m]}_t}}&\leq\Dt_1(p) \brak{\Ebh_t\edg{\exp\crl{3p\gamma\sigmat^2\brak{\abv{\xi}+\int_t^T\brak{\ha_s+\frac\gamma2+\lm\abv{Y^{[m-1]}_s}}ds}}}}^n\\
	&\quad\x\prod_{j=1}^n\exp\crl{3p\gamma\sigmat^2\Cs\sup_{s\in[t,T]}\Ebh\edg{\abv{\Yb^{[m],j}_s}}}.
\end{align*}
From Corollary \ref{corollary:Doob} and Jensen's inequality, we have
\begin{align*}
	&\quad\Ebh\edg{\exp\crl{3p\gamma\sigmat^2\sup_{s\in[t,T]}\abv{Y^{[m]}_s}}}\\
	&\leq \Ah_G	\Dt_1(p)\Ebh\edg{\exp\crl{6n(1+\Cs)p\gamma\sigmat^2\brak{\abv{\xi}+\int_t^T\brak{\ha_s+\frac\gamma2+\lm\abv{Y^{[m-1]}_s}}ds}}}\\
	&\leq \brak{\Ds(\xi,p)\Ebh\edg{\exp\crl{12n(1+\Cs)p\gamma\sigmat^2\lm(T-t)\sup_{s\in[t,T]}\abv{Y^{[m-1]}_s}}}}^\half,
\end{align*}
where 
\begin{align*}
	&\quad \Ds(\xi,p)\\
	&:=\brak{\Ah_G\Dt_1(p)}^2\brak{\Ebh\edg{\exp\crl{24n(1+\Cs)p\gamma\sigmat^2\abv{\xi}}}}^\half\brak{\Ebh\edg{\exp\crl{24n(1+\Cs)p\gamma\sigmat^2\brak{\int_0^T\ha_sds+\frac{\gamma T}2}}}}^\half.
\end{align*}
Choosing $h_5\in(0,T]$ such that $4n(1+\Cs)\lm h_5<1$ and $T=N_h^5h_5$ for some $N_h^5\in\N^*$, we get
\begin{align}
	&\Ebh\edg{\exp\crl{3p\gamma\sigmat^2\sup_{s\in[T-h_5,T]}\abv{Y^{[m]}_s}}}\label{eq:Ymd1}\\
	&~~~~~~~~~~~~~~~~~~~\leq \brak{\Ds(\xi,p)}^{\half+\frac14+\cdots+\frac1{2^m}}\brak{\Ebh\edg{\exp\crl{3p\gamma\sigmat^2\sup_{s\in[T-h_5,T]}\abv{Y^{[0]}_s}}}}^{\frac1{2^m}}\leq \Ds(\xi,p).\notag
\end{align}
Iterating backwardly, we get on $[(k-1)h_5,kh_5],~k=N_h^5-1,\cdots,1$,
\begin{align}
	&\quad\Ebh\edg{\exp\crl{3p\gamma\sigmat^2\sup_{s\in[(k-1)h_5,kh_5]}\abv{Y^{[m]}_s}}}\label{eq:Ymd2}\\
	&\leq \brak{\Ah_G\Dt_1(p)}^2\brak{\Ebh\edg{\exp\crl{24n(1+\Cs)p\gamma\sigmat^2\abv{Y^{[m]}_{kh_5}}}}}^\half\brak{\Ebh\edg{\exp\crl{24n(1+\Cs)p\gamma\sigmat^2\brak{\int_0^T\ha_sds+\frac{\gamma T}2}}}}^\half\notag\\
	&\leq\brak{\Ah_G\Dt_1(p)}^2\brak{\Ds(Y^{[m]}_{(k+1)h_5},\lceil8n(1+\Cs)\rceil p)}^\half\brak{\Ebh\edg{\exp\crl{24n(1+\Cs)p\gamma\sigmat^2\brak{\int_0^T\ha_sds+\frac{\gamma T}2}}}}^\half\notag\\
	&=:\Ds(Y^{[m]}_{kh_4},p),\notag
\end{align}
and derive a uniform estimate, i.e. $$\sup_{m\in\N^*}\Ebh\edg{\exp\crl{3p\gamma\sigmat^2\sup_{s\in[0,T]}\abv{Y^{[m]}_s}}}\leq\Dsh(\xi,p),$$ analogously to $\Dch(\xi,p)$ in Step 1 of Theorem \ref{theorem:main}.

\noindent{\bf Convergence.}
In order to prove that $\brak{Y^{[m]}}_{m\in\N^*}$ is a Cauchy sequence, we apply a $\theta$-method as Steps 2-3 of Theorem \ref{theorem:main}. Fix $m,k\geq1$ and $\theta\in(0,1)$, we define for $\phi=Y,~\Yb,~R^T-R,~\xi$,
\begin{align*}
	\dt\phi^{[m,q]}:=\frac{\phi^{[m+q]}-\theta\phi^{[m]}}{1-\theta},~\dtt\phi^{[m,q]}:=\frac{\phi^{[m]}-\theta\phi^{[m+q]}}{1-\theta},~\mbox{and}~~\dta\phi^{[m,q]}:=\abv{\dt\phi^{[m,q]}}+\abv{\dtt\phi^{[m,q]}},
\end{align*}
and for each $1\leq j\leq n$,
\begin{align*}
	\dt\phi^{[m,q],j}:=\frac{\phi^{[m+q],j}-\theta\phi^{[m],j}}{1-\theta},~\dtt\phi^{[m,q],j}:=\frac{\phi^{[m],j}-\theta\phi^{[m+q],j}}{1-\theta},~\mbox{and}~~\dta\phi^{[m,q],j}:=\abv{\dt\phi^{[m,q],j}}+\abv{\dtt\phi^{[m,q],j}}.
\end{align*}
Analogously, we can select $h_6\in(0,T]$ such that $64n(1+2\Cs)\lm h_6<1$ and $T=N_h^6h_6$ for some $N_h^6\in\N^*$. So we have
\begin{align*}
	&\quad\Ebh\edg{\exp\crl{3p\gamma\sigmat^2\sup_{s\in[T-h_6,T]}\dta Y^{[m,q]}_s}}\leq \brak{\Ab_G^\dagger(\xi,p)}^{\frac14}\brak{\Ebh\edg{\exp\crl{3p\gamma\sigmat^2\sup_{s\in[T-h_6,T]}\dta Y^{[m-1,q]}_s}}}^{\frac14}\\
	&\leq \brak{B_G^*(\xi,p)}^{\frac14+\cdots+\frac1{4^{m-1}}}\brak{\Ebh\edg{\exp\crl{3p\gamma\sigmat^2\sup_{s\in[T-h_6,T]}\dta Y^{[1,q]}_s}}}^{\frac1{4^{m-1}}}\\
	&\leq \brak{B_G^*(\xi,p)}^{\frac14+\cdots+\frac1{4^{m}}}\brak{\Ebh\edg{\exp\crl{\frac{6}{1-\theta}p\gamma\sigmat^2\sup_{s\in[T-h_6,T]} \abv{Y^{[q]}_s}}}}^{\frac1{4^{m}}},
\end{align*}
where $B_G^*(\xi,p)$ is dependent on $(n,\hv,\lv,\uv)$ via Lemma \ref{lemma:Rtheta} and bounded in analogy to $A^*_G(\xi,p)$.
Then we can derive a local estimate, i.e. for any $\theta\in(0,1)$, $$\lim_{m\to\infty}\sup_{q\in\N^*}\Ebh\edg{\exp\crl{3p\gamma\sigmat^2\sup_{t\in[T-h_6,T]}\dta Y^{[m,q]}_s}}\leq B_G^*(\xi,p).$$ In the thoughts of
estimates \eqref{eq:Ymd1}-\eqref{eq:Ymd2} and local results on small intervals following Theorem 2.8 in \cite{MDQBSDE}, we can get the desired convergence and verification similarly to Steps 2-3-4 in Theorem \ref{theorem:main}. Moreover, the uniqueness can be obtained identically to that in \cite{MDQBSDE,LiumdGBSDE}.
\qed

\appendix
\section{Appendix}
This section is dedicated to some a priori estimates of standard $G$-BSDEs, essential to the establishments of DMR$G$-BSDEs in Section 4. We consider $G$-BSDEs taking the form:
\begin{equation}\label{eq:GBSDEappendix}
	Y_t=\xi+\int_t^Tf(s,Y_s,Z_s)ds-\int_t^TZ_sdB_s-(K_T-K_t),~t\in[0,T],
\end{equation}
\subsection{Lipschitz $G$-BSDEs}

\begin{Proposition}[\cite{HuJiPengSong}]\label{prop:1daprDY}
	Let $(\xi^{(m)}, f^{(m)},g_{ij}^{(m)}),~ m=1,2$ satisfy Assumption \ref{assumption:HLip}. Suppose that $(Y^{(m)},Z^{(m)},K^{(m)})$ is a solution to $G$-BSDE \eqref{eq:GBSDEappendix} with data $(\xi^{(m)},f^{(m)},g_{ij}^{(m)})$, respectively. Set $\Yh_t:=Y^{(1)}_t-Y^{(2)}_t$. Then for $1< \a\leq\b$, there is some constant $C>0$ depending on $T$, $\sigmat$, $\sigmah$, $L$ and $\a$ such that
	$$\abv{\Yh_t}^\a\leq C\Ebh_t\edg{\abv{\hat{\xi}}^\a+\brak{\int_t^T\hat{h}_sds}^\a},$$	
	where $\hat{\xi}=\xi^{(1)}-\xi^{(2)}$ and
	
	$\hat{h}_s=\abv{f^{(1)}(s,Y_s^{(2)},Z_s^{(2)})-f^{(2)}(s,Y_s^{(2)},Z_s^{(2)})}+\sum_{i,j=1}^d\abv{g_{ij}^{(1)}(s,Y_s^{(2)},Z_s^{(2)})-g_{ij}^{(2)}(s,Y_s^{(2)},Z_s^{(2)})}$.
\end{Proposition}
\begin{Proposition}[\cite{HuJiPengSong}]\label{prop:1daprDZ}
	Let $(\xi^{(m)}, f^{(m)},g_{ij}^{(m)}),~ m=1,2$ satisfy Assumption \ref{assumption:HLip} for the same $\b>1$ and $L>0$. Suppose that $(Y^{(m)},Z^{(m)},K^{(m)})\in \Sc_G^\a(0,T)$, for some $1\leq\a<\b$, is a solution to $G$-BSDE \eqref{eq:GBSDEappendix} with data $(\xi^{(m)},f^{(m)},g_{ij}^{(m)})$, respectively. Set $\Yh_t:=Y^{(1)}_t-Y^{(2)}_t$ and $\Zh_t:=Z^{(1)}_t-Z^{(2)}_t$. Then there exists $C>0$ depending on $T$, $\sigmat$, $\sigmah$, $L$ and $\a$ such that 
	$$\Ebh\edg{\brak{\int_0^T\abv{\Zh_s}^2ds}^\ah}\leq C\crl{\norm{\Yh}_{S_G^\a}^\a+\norm{\Yh}_{S_G^\a}^\ah\sum_{m=1}^2\edg{\norm{Y^{(m)}}_{S_G^\a}^\ah+\norm{\int_0^Th_s^{(m)}ds}_{L_G^\a(\Om_T)}^\ah}},$$
	where $h_s^{(m)}=\abv{f^{(m)}(s,0,0)}+\sum_{i,j=1}^d\abv{g_{ij}^{(m)}(s,0,0)} $.
\end{Proposition}

\subsection{Quadratic $G$-BSDEs}
In this section, we list some useful a priori estimates on Q$G$-BSDEs, cf. \cite{QGBSDEbt,QGBSDEubt,2QGBSDE} and so on.
The unique existence under Assumptions \ref{assumption:HQbdd} of $\Sc_G^2(0,T)$-solution $(Y,Z,K)$ to Q$G$-BSDE \eqref{eq:GBSDEappendix} can be derived by Theorem 5.3 in \cite{QGBSDEbt}. From Lemma 3.1 and Theorem
3.2 in \cite{2QGBSDE}, there exist some upper bounds such that
\begin{equation}\label{eq:qgbsdeyzapr}
	\norm{Y}_{S_G^\infty(0,T)}+\norm{Z}_{\mathrm{BMO}_G}\leq C_{YZ}:=C_{YZ}(T,M_0,\lm,\gamma,\sigmat,\sigmah),
\end{equation}
and 
\begin{equation}\label{eq:qgbsdekapr}
	\Ebh\edg{\abv{K_T}^p}\leq C_K:=C_K(T,p,M_0,\lm,\gamma,\sigmat,\sigmah),~p\geq1,
\end{equation}
which implies that $(Y,Z,K)\in S_G^\infty(0,T)\x{\rm BMO}_G\x\As_G^T$.

\begin{Proposition}\label{proposition:qgbsdeaprdyz}
	Let $(\xi^i,f^i),~i=1,2$ satisfy Assumption \ref{assumption:HQbdd} for the same positive constants $M_0, ~\lm$ and $\gamma$. Suppose that $(Y^i,Z^i,K^i)\in S_G^2(0,T)\x H_G^2(0,T;\R^d)\x\Ac_G^2(0,T)$ are solutions to Q$G$-BSDEs \eqref{eq:GBSDEappendix} with data $(\xi^i,f^i)$ respectively. Then for any $t\in[0,T]$,
	\begin{equation*}
		\Ebh\edg{\brak{\int_t^T\abv{Z_s^1-Z^2_s}^2ds}^{\frac p2}}\leq C_{\Zh} \crl{\norm{\xi^1-\xi^2}_{L_G^\infty}^p+\Ebh\edg{\sup_{s\in[t,T]}\abv{Y^1_s-Y^2_s}^p}},~p\geq2,
	\end{equation*}
	where $C_{\Zh}=C_{\Zh}(T,p,\sigmat,\sigmah,M_0,\lm,\gamma)>0$.
\end{Proposition}
\noindent{\bf\underline{Proof.}}
We apply It\^o's formula to $\abv{Y^1-Y^2}^2$ and $\abv{Y^i}^2$, we have by estimates \eqref{eq:qgbsdeyzapr}-\eqref{eq:qgbsdekapr} and Yong's inequality, we can get the desired result.
\qed
\begin{Remark}\label{remark:QGZ}
	In view of estimate \eqref{eq:qgbsdeyzapr} and Proposition \ref{proposition:qgbsdeaprdyz}, if $(Y,Z,K)\in S_G^2(0,T)\x H_G^2(0,T;\R^d)\x\Ac_G^2(0,T)$ is a solution to Q$G$-BSDE \eqref{eq:GBSDEappendix} under Assumption \ref{assumption:HQbdd}, then
	\begin{equation*}
		(Y,Z,K)\in S_G^\infty(0,T)\x\Hs_G^T(\R^d)\x\As_G^T.
	\end{equation*}
\end{Remark}
\begin{Lemma}[\cite{QGBSDEubt}]\label{lemma:QGBSDEaprY}
	Assume that $(Y,Z,K)\in S_G^2(0,T)\x H_G^2(0,T;\R^d)\x\Ac_G^2(0,T)$ satisfies the following equation
	\begin{equation*}
		Y_t=\xi+(\Kt_T-\Kt_t)+\int_t^Tf(s,Z_s)ds-\int_t^TZ_sdB_s-(K_T-K_t),~t\in[0,T],
	\end{equation*}
	with $\Kt\in \Ac_G^2(0,T)$. Suppose that there are two constants $p\geq1$ and $\eps>0$ such that
	\begin{equation}
		\Ebh\edg{\exp\crl{(2p+\eps)\kappa\sigmat^2\sup_{t\in[0,T]}\abv{Y_t}+(2p+\eps)\kappa\sigmat^2\int_0^T\hb_tdt}}<\infty.
	\end{equation}
	\begin{enumerate}[($i$)]
		\item If $\abv{f(t,\om,z)}\leq\hb_t(\om)+\frac\kappa2\abv{z}^2$, then for each $t\in[0,T]$,
		\begin{equation*}
			\exp\crl{p\kappa\sigmat^2\abv{Y_t}}\leq\Ebh_t\edg{\exp\crl{p\kappa\sigmat^2\brak{\abv{\xi}+\int_t^T\hb_sds}}}.
		\end{equation*}
		\item If $f(t,\om,z)\leq\hb_t(\om)+\frac{\kappa}{2}\abv{z}^2$, then for each $t\in[0,T]$,
		\begin{equation*}
			\exp\crl{p\kappa\sigmat^2Y_t^+}\leq\Ebh_t\edg{\exp\crl{p\kappa\sigmat^2\brak{\xi^++\int_t^T\hb_sds}}}.
		\end{equation*}
	\end{enumerate} 
\end{Lemma}

\end{document}